\documentclass{article}
\usepackage{graphicx} 
\usepackage{amsmath}
\usepackage{algorithm}
\usepackage{algpseudocode}
\usepackage{amssymb}
\usepackage{fullpage}
\usepackage{xcolor}
\usepackage{bm}
\usepackage{booktabs}
\usepackage{array}
\usepackage{amsthm}
\usepackage[numbers,sort&compress]{natbib}
\usepackage{url}
\usepackage{placeins}
\usepackage{etoolbox}
\AtBeginEnvironment{table}{\setlength{\abovecaptionskip}{4pt}\setlength{\belowcaptionskip}{6pt}}

\newtheorem{proposition}{Proposition}

\title{A Compact Residual-Screened Local Reduced-Order Model\\
for Diffuse-Domain Shape Optimization}
\author{Xiaoping Zhang\thanks{School of Mathematics and Statistics, Wuhan University, Wuhan, Hubei, 430072, P.R. China.} \and Anusooya Vasudevan\thanks{Department of Mathematics, University of South Carolina, Columbia, SC 29208, USA.} \and Lili Ju\footnotemark[2] \and Zhu Wang\footnotemark[2]\ \thanks{Corresponding author. Email: \texttt{wangzhu@math.sc.edu}}}
\date{}

\begin{document}

\maketitle

\begin{abstract}
We generalize the compact residual-screened local reduced-order model (CL-ROM) from Brinkman-penalized flow to shape optimization with the diffuse domain method (DDM). As with Brinkman penalization, DDM extends evolving-domain problems to a fixed rectangular domain, allowing one grid for discretization and directly compatible snapshots from different designs. Unlike Brinkman penalization, DDM uses a phase field to weight operators and sources and treats Dirichlet, Neumann and Robin conditions through diffuse boundary terms. The CL-ROM mechanism carries over: recent full-order and corrected solutions determine local reduced spaces, while residual screening, correction with retained factors and objective verification decide when a full-order solve is needed. We derive reduced operators, residuals and sensitivities for each diffuse model, using Galerkin projection for the symmetric problems and least-squares Petrov--Galerkin projection for nonlinear flow. For the symmetric problems we prove a local approximation estimate that admits inexact, corrected snapshots. In numerical tests of shape identification, drag minimization in Navier--Stokes channel flow and compliance minimization of a cantilever, CL-ROM uses reduced solutions for nearly all design updates, reproduces the full-order objectives with relative differences below $5\times10^{-6}$, and is $3.85$ to $7.65$ times faster than a full-order baseline that reuses factorizations as preconditioners.
\end{abstract}

\noindent\textbf{Keywords:} diffuse domain method; shape optimization;
phase field; local reduced-order model
\vspace{2mm}

\section{Introduction}
\label{sec:intro}

In shape optimization, one seeks the geometry of a domain, or of a body or structure inside it, that minimizes a cost functional subject to partial differential equations (PDEs) posed on that geometry \cite{pironneau1984optimal,sokolowski1992introduction}. It has broad applications, including the design of bodies with low drag in viscous flow \cite{mohammadi2010applied,gunzburger2002perspectives,borrvall2003topology}, the layout of stiff and light structures \cite{bendsoe1988generating,bendsoe2003topology,allaire2004structural,takezawa2010shape}, and the identification of an unknown domain from observed data \cite{laporte2003numerical,burger2005survey}. Most methods are gradient based and update the design with shape sensitivities computed from a state and an adjoint solution \cite{sokolowski1992introduction,giles2000introduction}, so each design iteration typically requires state and adjoint PDE solves on the current geometry. Since an optimization may take hundreds or thousands of iterations, these solves dominate the total cost. With body-fitted discretizations, a change of geometry may in addition require remeshing and the transfer of solution fields between meshes, which further increases the cost.

Reduced-order models (ROMs) provide computationally efficient surrogates for parametric PDEs \cite{benner2015survey,benner2017model}. A ROM is typically built by first extracting a low-dimensional reduced basis from high-fidelity full-order solution snapshots, for example by proper orthogonal decomposition (POD) \cite{holmes2012turbulence,kunisch2001galerkin,iliescu2014are} or reduced-basis methods \cite{hesthaven2016certified,quarteroni2015reduced}, and then seeking approximate solutions in the span of that basis, either by projection (intrusive ROMs) \cite{benner2015survey,carlberg2011efficient} or by inference from data (non-intrusive ROMs) \cite{kramer2024learning}. The resulting ROM can be evaluated rapidly at new parameter values. Shape optimization, however, is a difficult setting for ROMs. First, snapshots computed on different geometries, and hence on different body-fitted meshes, do not belong to one discrete vector space and cannot be combined directly into a basis. Second, solution features such as boundary layers and wakes move with the boundary; such transport-like dependence on the design is poorly captured by a single linear space of low dimension \cite{peherstorfer2022breaking}. Third, a parameterized design may have many degrees of freedom and its trajectory is not known in advance, so an offline training set that covers the designs visited by the optimizer can be prohibitively large.

Several strategies address these difficulties. Mapping each geometry to a reference domain places all snapshots on a common mesh, at the cost of geometry-dependent coefficients in the transformed operators and a parametrization of the geometry \cite{quarteroni2015reduced}; combined with free-form deformation, this approach has been used for reduced-basis shape optimization \cite{manzoni2012shape}. Embedded discretizations instead keep a fixed background mesh; reduced models have been built for embedded boundary models with evolving interfaces \cite{balajewicz2014reduction}, the shifted boundary method \cite{karatzas2020reduced}, and a diffuse-interface Stokes--Darcy model \cite{stoter2022deim}, where they are combined with discrete empirical interpolation \cite{chaturantabut2010deim}. Local reduced bases restrict each basis to a part of the solution manifold \cite{amsallem2012nonlinear}, and adaptive ROMs build or update the reduced space along the optimization path rather than from an offline training set, in progressive and trust-region frameworks \cite{zahr2015progressive,qian2017certified} and in topology optimization \cite{gogu2015improving,xiao2024primaldual}.

The companion work \cite{li2026residual} combined these last two ideas in the \emph{compact residual-screened local ROM} (CL-ROM) for Brinkman-penalized shape optimization of incompressible flows. Brinkman penalization embeds the solid body in a fixed rectangular domain through a design-dependent friction term, so all states share one structured grid and snapshots from different designs are directly compatible. CL-ROM builds local reduced state and adjoint spaces from short windows of recent full-order and corrected solutions. The full-order residual of a reduced solution decides whether it is accepted; a rejected solution is corrected with a retained factorization or replaced by a full-order solve, and objective verification, triggered by a dual-weighted error estimate, protects the objective. In that setting, however, the design enters the equations only through a diagonal friction coefficient, and the only boundary condition imposed on the body is no slip.

In this work, we extend CL-ROM to shape optimization based on the diffuse domain method (DDM). Like Brinkman penalization, DDM extends a PDE posed on a complex or evolving domain to a fixed rectangular domain by representing the domain with a phase field and its boundary with a diffuse interface \cite{li2009solving}. Boundary conditions, however, enter through problem-dependent volume weights, diffuse surface terms and penalties, so that Dirichlet, Neumann and Robin conditions are handled in one framework \cite{li2009solving,lervag2015analysis,burger2017analysis,yu2012extended}; in shape and topology optimization, phase-field representations are also widely used \cite{takezawa2010shape,garcke2015numerical}. We consider three representative problems: shape identification with Dirichlet data imposed by a variational surface penalty, drag minimization of a body in Navier--Stokes channel flow, with no slip imposed by a weighted volume penalty following \cite{vasudevan2026diffuse}, and compliance minimization of a cantilever whose traction-free boundary is represented by a phase-field-weighted stiffness \cite{yu2012extended,takezawa2010shape}. Together they cover prescribed-value and natural boundary conditions, linear and nonlinear state equations, and problems that require a separate adjoint solve as well as a self-adjoint problem that requires only a state solve.

The mechanism of CL-ROM carries over naturally. All states are again discretized on one fixed grid, so snapshots from different designs need neither extension nor mesh transfer, and the design changes slowly over long phases of the optimization, so recent snapshots represent the next state well. What changes is how the design enters the discrete equations. Beyond a diagonal penalty, which the flow model retains, the phase field enters through the weighted coefficients of the principal part of the operator, the diffuse surface measure of a boundary penalty, the load, and the nonlinear convection. The reduced operators, the full-order residuals and the sensitivities must account for each of these terms, and the projection used to build the ROM and its correction must respect the structure of each state system, which is symmetric positive definite for shape identification and compliance and a nonlinear saddle-point system for the flow.

The main contributions of this work are threefold. (1) We extend CL-ROM to DDM shape optimization, deriving the reduced operators, residuals and sensitivities of the three diffuse models, with Galerkin projection for the symmetric problems and least-squares Petrov--Galerkin projection for the nonlinear flow problem. (2) We derive a local approximation estimate for the symmetric problems that allows inexact, corrected snapshots and separates the effects of design change, snapshot error and POD truncation. (3) We show on the three problems that CL-ROM reproduces the full-order objectives and designs at a fraction of the cost, and we measure the errors of residual-screened candidates at fixed designs, finding that sampled objective errors remain below the prescribed tolerance even when calibrated residual indicators underestimate state errors.

The remainder of the paper is organized as follows. Section~\ref{sec:formulation} presents the three DDM shape-optimization problems, their optimality systems and discretizations, and the common design update. Section~\ref{sec:trrom} describes the extension of CL-ROM, its local approximation estimate and the problem-dependent ingredients. Section~\ref{sec:experiments} reports the numerical experiments on the three problems and the accuracy of accepted reduced solutions. Section~\ref{sec:discussion} concludes the paper and outlines directions for further reducing the cost.

\FloatBarrier
\section{Shape Optimization in the Diffuse-Domain Formulation}
\label{sec:formulation}

Three shape-optimization problems are considered in this work. For each,
we first use the diffuse-domain method to extend the state equation from the
evolving domain $\Omega$ to a fixed rectangular domain $D$, and then briefly
present the corresponding optimality system and discretization.
Section~\ref{sec:ddm} sets the common phase-field notation and uses the
Poisson equation to illustrate how Dirichlet, Neumann and Robin boundary
conditions enter DDM. Sections~\ref{sec:opt}--\ref{sec:elast}
then treat shape identification, drag minimization in a channel and
compliance minimization of a cantilever in turn. All three designs are
evolved by one gated Allen--Cahn flow, stated in Section~\ref{sec:design}
together with the stopping rule. The derivations are collected in the
appendices.

\subsection{Diffuse-domain constructions}
\label{sec:ddm}

\begin{figure}[htbp]
\centering
\includegraphics[width=0.88\textwidth]{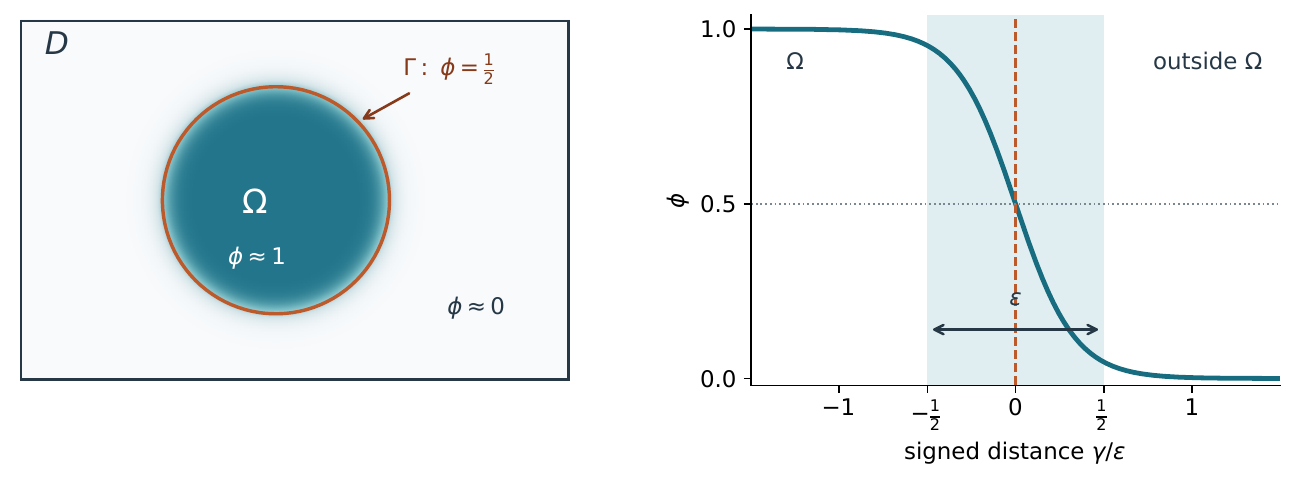}
\caption{Diffuse-domain geometry and the phase-field profile \eqref{eq:phase}.
The fixed box $D$ contains the design domain $\Omega$; its free boundary
$\Gamma$ is represented by the $\phi=1/2$ contour. Across the diffuse
interface, $\phi$ changes from approximately one to approximately zero
over a distance set by $\varepsilon$. The shaded interval in the profile
is $-\varepsilon/2\leq\gamma\leq\varepsilon/2$.}
\label{fig:ddm_schematic}
\end{figure}

Let $D \subset \mathbb{R}^d$, $d = 2, 3$, be a fixed box containing the design domain $\Omega$ with boundary $\Gamma=\partial\Omega$. In the shape-optimization problems below, this free boundary is denoted by $\Gamma_f$. In the DDM formulation, this geometry is encoded by a phase field $\phi : D \to [0,1]$ with the equilibrium profile
\begin{equation}
\label{eq:phase}
\phi(\bm{x}) = \frac{1}{2} \Big( 1 - \tanh \big( \frac{3\, \gamma(\bm{x})}{\varepsilon} \big) \Big),
\end{equation}
where $\gamma$ is the signed distance to $\Gamma$, negative in $\Omega$, and $\varepsilon$ is the interface width: $\phi \approx 1$ in $\Omega$ and $\phi \approx 0$ outside it, with a transition layer of width $O(\varepsilon)$. Figure~\ref{fig:ddm_schematic} illustrates this ideal profile. The surface measure of $\Gamma$ is approximated by the diffuse delta function $\delta_\varepsilon = |\nabla\phi|$, $\int_\Gamma v\,\mathrm{d}s \approx \int_D \delta_\varepsilon\,v\,\mathrm{d}\bm{x}$ for smooth $v$.

Consider, for example, the Poisson equation $-\Delta u=f$ in $\Omega$ with a
Dirichlet, Neumann or Robin condition on $\Gamma$, where $\bm{n}$ denotes
the unit outward normal of $\Omega$. DDM extends $f$ from $\Omega$ to $D$
and rewrites the weak form on $D$, replacing integrals over $\Omega$ by
$\phi$-weighted integrals and integrals over $\Gamma$ by
$\delta_\varepsilon$-weighted integrals. The resulting models are
\cite{li2009solving,burger2017analysis}:
\begin{equation}
\label{eq:ddm_bc}
\begin{aligned}
&\text{Dirichlet, } u = g_D: && -\nabla\cdot(\phi\nabla u) + \beta_\varepsilon^{-1}\delta_\varepsilon\,(u-g_D) = \phi f,\\
&\text{Neumann, } \partial_{\bm{n}} u = g_N: && -\nabla\cdot(\phi\,\nabla u) = \phi f + \delta_\varepsilon\, g_N,\\
&\text{Robin, } \partial_{\bm{n}} u + \kappa\,(u - u_\infty) = 0: && -\nabla\cdot(\phi\,\nabla u) + \kappa\,\delta_\varepsilon\,(u - u_\infty) = \phi f,
\end{aligned}
\end{equation}
where $\kappa>0$ and $u_\infty$ are the transfer coefficient and ambient
value, and $\beta_\varepsilon>0$ is a small penalty parameter. In the
interior, where $\phi\approx1$ and $\delta_\varepsilon\approx0$, each
model reduces to the Poisson equation, so the boundary condition acts only
through the $\delta_\varepsilon$-weighted terms near $\Gamma$. Neumann and
Robin data enter through their physical flux terms and need no artificial
parameter. The Dirichlet model is the Robin model with
$\kappa=\beta_\varepsilon^{-1}$ and $u_\infty=g_D$: the condition is
imposed by the penalty $\beta_\varepsilon^{-1}\delta_\varepsilon(u-g_D)$,
which enforces $u\approx g_D$ on $\Gamma$ for small $\beta_\varepsilon$.

In the shape optimization problems, the phase field becomes the design
variable, so changes of shape, including splitting or merging of
components, require no remeshing. Equation~\eqref{eq:phase} then
describes an ideal profile, not a constraint on the evolving design: we
continue to define $\Gamma$ as the $\phi=1/2$ contour, although the
transition need not remain a signed-distance profile.

\begin{figure}[htbp]
\centering
\includegraphics[width=\textwidth]{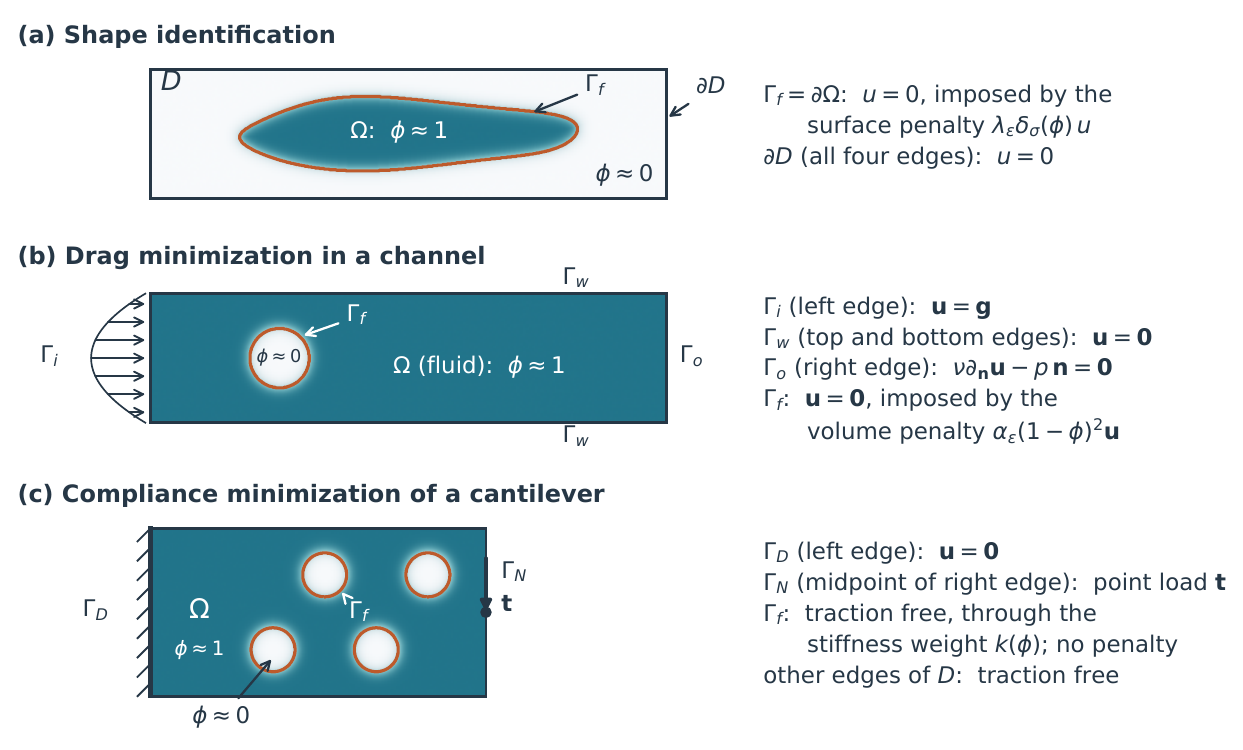}
\caption{Boundary conditions and phase-field conventions of the three DDM
state equations, drawn for the configurations of
Section~\ref{sec:experiments} (not to scale). Shading shows $\phi$, from
$\phi\approx1$ in $\Omega$ (dark) to $\phi\approx0$ outside it (light); the
red curve is the $\phi=\tfrac12$ contour, which represents the free boundary
$\Gamma_f$, the part of $\partial\Omega$ that moves with the design. The
fixed boundary pieces lie on $\partial D$. The lists give the condition on
each piece and, on $\Gamma_f$, the DDM term that imposes it.
(a)~Shape identification (Sections~\ref{sec:opt} and \ref{sec:num}): all of $\partial\Omega$
is free. (b)~Drag minimization in a channel (Sections~\ref{sec:flows} and \ref{sec:flowexp}):
$\Omega$ is the fluid and $\Gamma_f$ the body surface. (c)~Compliance
minimization of a cantilever (Sections~\ref{sec:elast} and \ref{sec:cantexp}): $\Omega$ is the
material and $\Gamma_f$ its traction-free boundary.}
\label{fig:settings}
\end{figure}

\subsection{Shape identification}
\label{sec:opt}

We first consider a shape-identification problem, which aims to determine the
boundary of a domain so that the solution of the Poisson equation on
it matches a given observation. In the
synthetic experiment, $\bar u=\chi_{\Omega_{\rm ref}}u_{\rm ref}$, where
$u_{\rm ref}$ solves the Poisson equation on the prescribed reference domain
$\Omega_{\rm ref}$; hence $\bar u=0$ outside $\Omega_{\rm ref}$. Using $\bar u$ as data, we seek a domain $\Omega$ whose entire boundary
$\Gamma_f=\partial\Omega$ is free (Figure~\ref{fig:settings}(a)) by solving
\begin{equation}
\label{eq:objective}
\begin{aligned}
\min_{\Omega\subset D}\quad &J(\Omega)=\int_D(\chi_{\Omega} u(\Omega)-\bar u)^2\,\mathrm{d}\bm{x},\\
\text{s.t.}\quad &-\Delta u=f\ \text{in }\Omega,\qquad u=0\ \text{on }\Gamma_f.
\end{aligned}
\end{equation}

\paragraph{Diffuse-domain formulation.}
The free boundary $\Gamma_f$ has homogeneous Dirichlet data, so we use the Dirichlet model of \eqref{eq:ddm_bc}, the variational surface penalty. We replace its principal-part weight $\phi$ by $k(\phi)=k_{\min}+(1-k_{\min})\phi$, whose small floor $k_{\min}>0$ keeps the discrete operator nonsingular, and its surface measure $\delta_\varepsilon$ by the regularized geometric measure
\begin{equation*}
\delta_\sigma(\phi)=\sqrt{|\nabla\phi|^2+\sigma^2}-\sigma,
\qquad
\lambda_\varepsilon=\beta_\varepsilon^{-1}=(c_D\varepsilon^{3/4})^{-1},
\end{equation*}
with a small regularization $\sigma>0$ and a constant $c_D>0$. The scaling
$\beta_\varepsilon=c_D\varepsilon^{3/4}$ is the one for which convergence as $\varepsilon\to0$ is
proved on a fixed smooth domain \cite{burger2017analysis}. The floor $k_{\min}$ and the regularization $\sigma$, which
keeps the derivative of the surface measure defined where $\nabla\phi=0$,
are small perturbations of the analyzed model; Appendix~\ref{app:ddm_refinement}
verifies that the implemented model converges at approximately second order
under mesh refinement at fixed interface width. We set
$u=0$ on the fixed boundary $\partial D$; the surface penalty imposes zero
Dirichlet data on the free boundary $\Gamma_f$
(Figure~\ref{fig:settings}(a)). The optimization problem is
\begin{equation*}
\min_\phi J(\phi)=\int_D(\phi u-\bar u)^2\,\mathrm{d}\bm{x}
\quad\text{subject to}\quad
-\nabla\cdot(k(\phi)\nabla u)+\lambda_\varepsilon\delta_\sigma(\phi)u=\phi f
\ \text{in }D,\qquad u=0\ \text{on }\partial D .
\end{equation*}
\paragraph{Optimality system.}
Let $w$ be the adjoint with zero Dirichlet boundary condition on $\partial D$, and let $e=\phi u-\bar u$. The state operator is symmetric, so the adjoint equation has the same operator with a different right-hand side:
\begin{equation}
\label{eq:adjoint_cont}
-\nabla\cdot(k(\phi)\nabla w)+\lambda_\varepsilon\delta_\sigma(\phi)w
=2\phi e .
\end{equation}
For a design variation $\vartheta$, the sensitivity is
\begin{equation}
\label{eq:sens_cont}
J'(\phi)[\vartheta]=\int_D
\big[2ue+fw-(1-k_{\min})\nabla u\cdot\nabla w\big]\vartheta
-\lambda_\varepsilon uw\,
\frac{\nabla\phi\cdot\nabla\vartheta}{\sqrt{|\nabla\phi|^2+\sigma^2}}
\,\mathrm{d}\bm{x}.
\end{equation}
The last term accounts for the motion of the diffuse surface measure. Appendix~\ref{app:derivation} derives the adjoint and sensitivity from the Lagrangian.

\paragraph{Discretization.}
We use a flux-form finite-difference scheme on a uniform Cartesian grid of spacing $h$, with state and design values stored at grid nodes. Let $G$ map nodal values to face differences, $\Pi_f$ average nodal coefficients to faces, and $D_x,D_y$ form centered nodal differences for the diffuse surface measure. These fixed sparse operators include the homogeneous outer boundary. With $\bm\phi$ the vector of interior nodal design values and $\delta_\sigma(\bm\phi)=\sqrt{(D_x\bm\phi)^2+(D_y\bm\phi)^2+\sigma^2}-\sigma$ entrywise, the discrete state system is
\begin{equation}
\label{eq:state_discrete}
A(\bm\phi)\mathbf U=\mathbf b(\bm\phi),\qquad
A(\bm\phi)=G^\top\operatorname{diag}\!\big(k_{\min}+(1-k_{\min})\Pi_f\bm\phi\big)G
+h^2\lambda_\varepsilon\operatorname{diag}(\delta_\sigma(\bm\phi)),
\qquad \mathbf b(\bm\phi)=h^2\bm\phi\circ\mathbf f.
\end{equation}
The objective is $J_h=h^2\|\bm\phi\circ\mathbf U-\bar{\mathbf U}\|^2$. Writing $\mathbf e=\bm\phi\circ\mathbf U-\bar{\mathbf U}$, the exact discrete adjoint is
\begin{equation*}
A(\bm\phi)\mathbf W=\nabla_{\mathbf U}J_h=2h^2\bm\phi\circ\mathbf e .
\end{equation*}
The gradient combines the explicit derivative of $J_h$ with the transposed design derivatives of \eqref{eq:state_discrete}, applied to $\mathbf U$ and $\mathbf W$:
\begin{equation}
\label{eq:sensitivity}
\begin{aligned}
\nabla_{\bm\phi}J_h={}&2h^2\mathbf U\circ\mathbf e+h^2\mathbf f\circ\mathbf W
-(1-k_{\min})\Pi_f^\top\!\big((G\mathbf U)\circ(G\mathbf W)\big)\\
&-h^2\lambda_\varepsilon\!\left[D_x^\top\!\left(\frac{\mathbf U\circ\mathbf W\circ D_x\bm\phi}
{\sqrt{(D_x\bm\phi)^2+(D_y\bm\phi)^2+\sigma^2}}\right)
+D_y^\top\!\left(\frac{\mathbf U\circ\mathbf W\circ D_y\bm\phi}
{\sqrt{(D_x\bm\phi)^2+(D_y\bm\phi)^2+\sigma^2}}\right)\right].
\end{aligned}
\end{equation}

\subsection{Drag minimization in a channel}
\label{sec:flows}

We next optimize the shape of a solid body immersed in a steady
incompressible flow. This problem has nonlinear convection and a coupled
velocity--pressure system. Let $\Omega\subset D$ be the fluid region
surrounding the body. Its boundary is
$\partial\Omega=\Gamma_i\cup\Gamma_o\cup\Gamma_w\cup\Gamma_f$, where
$\Gamma_i$ is the inflow, $\Gamma_o$ the outflow, $\Gamma_w$ the fixed
channel walls and $\Gamma_f$ the body surface, which is the free design
boundary (Figure~\ref{fig:settings}(b)). For a given body shape, the velocity $\bm{u}$ and pressure $p$ of
a Newtonian fluid with viscosity $\nu$ satisfy the steady incompressible
Navier--Stokes equations \cite{layton2008introduction}
\begin{equation}
\label{eq:nse_sharp}
\begin{gathered}
\begin{aligned}
-\nu \Delta \bm{u}+(\bm{u}\cdot\nabla)\bm{u}+\nabla p
  &=\bm{0} && \text{in }\Omega,\\
\nabla\cdot\bm{u} &=0 && \text{in }\Omega.
\end{aligned}
\\[-1mm]
\bm{u}=\bm{g}\ \text{on }\Gamma_i,\qquad
\bm{u}=\bm{0}\ \text{on }\Gamma_w\cup\Gamma_f,\qquad
\nu\,\partial_{\bm n}\bm{u}-p\bm n=\bm{0}\ \text{on }\Gamma_o.
\end{gathered}
\end{equation}
We seek the body boundary $\Gamma_f$ that minimizes the drag on the body
subject to a prescribed body volume, equivalently a prescribed fluid volume
in the fixed channel.

\paragraph{Diffuse-domain formulation.}
We use the weighted diffuse-domain flow construction of \cite{vasudevan2026diffuse}, with the modified no-slip penalty specified below. We represent $\Omega$ by the phase field ($\phi \approx 1$ in the fluid, $\phi \approx 0$ in the solid) and solve the following equations on $D$:
\begin{equation}
\label{eq:nse_ddm}
\phi\,\big( -\nu \Delta \bm{u} + \nabla p \big)
 + (\bm{u}\cdot\nabla)(\phi\,\bm{u})
 + \alpha_\varepsilon(1-\phi)^2\bm{u} = \bm{0},
\qquad \nabla\cdot(\phi\,\bm{u}) = 0 \qquad \text{in } D,
\end{equation}
where $\alpha_\varepsilon=C_p\nu/\varepsilon^2$ with a penalty constant $C_p>0$. The inflow, wall and outflow conditions of \eqref{eq:nse_sharp} are imposed on the corresponding fixed portions of $\partial D$ (Figure~\ref{fig:settings}(b)). The penalty damps velocity in the solid; squaring $1-\phi$ reduces its tail on the fluid side of the interface. The weighted continuity equation conserves the diffuse fluid flux $\phi\bm{u}$. The design enters the momentum, penalty, continuity and convection terms.

The choice in \eqref{eq:nse_ddm} to place $\phi$ outside the viscous and
pressure-gradient derivatives and impose no slip through a volume penalty is
not forced by the boundary condition. A divergence-form viscous operator with
a suitable Dirichlet penalty is also possible, but
$-\nu\nabla\cdot(\phi\nabla\bm{u})=-\nu\phi\Delta\bm{u}
-\nu(\nabla\phi\cdot\nabla)\bm{u}$ contains an additional interface term;
writing the pressure term as $\nabla(\phi p)$ likewise adds $p\nabla\phi$.
Such a formulation would change the state equation, its adjoint and the
design sensitivity. We retain \eqref{eq:nse_ddm} to test the ROM on a
nonlinear saddle-point system whose momentum and continuity equations both
depend on the design; we do not claim that this model is superior to a
divergence-form alternative. It also differs from the companion Brinkman
embedding \cite{li2026residual}, which leaves the differential operators
unweighted and introduces the design through a friction term.

\paragraph{Optimality system.}
The weighted momentum equation is not in conservation form, so its
penalty reaction alone does not represent body drag. Instead, we evaluate
the complete momentum flux through the fixed channel boundary:
\begin{equation}
\label{eq:dragflow}
J(\phi)
= \int_{\Gamma_i} \big(p+u_x^2-\nu\partial_xu_x\big)\,\mathrm{d}s
 -\int_{\Gamma_o} \big(p+u_x^2-\nu\partial_xu_x\big)\,\mathrm{d}s
 +\nu\int_{\Gamma_w}\partial_{\bm n}u_x\,\mathrm{d}s .
\end{equation}
For the sharp problem this equals body drag by momentum conservation;
evaluated on a diffuse solution, it is the drag approximation used here.
Here $u_x$ is the streamwise velocity component and the normal $\bm n$ points out of the channel. Wall friction therefore
subtracts from the inferred body drag, since $\partial_{\bm n}u_x<0$
on these walls. At the outlet, $p-\nu\partial_xu_x=0$ by the imposed
traction condition. A prescribed body volume is imposed through the
equivalent fluid-volume constraint
\begin{equation}
\label{eq:volcon}
\int_D\phi\,\mathrm{d}\bm{x}=V_0 .
\end{equation}
We minimize \eqref{eq:dragflow} subject to \eqref{eq:nse_ddm} and
\eqref{eq:volcon}. The design is fixed at $\phi=1$ near $\partial D$.
For smooth fields, the continuous adjoint $(\bm w,\pi)$ satisfies
\begin{equation}
\label{eq:adjoint_flow}
\begin{aligned}
-\nu\Delta(\phi\bm w)+[\nabla(\phi\bm u)]^{\!\top}\bm w
 -\phi\big((\bm u\cdot\nabla)\bm w+(\nabla\cdot\bm u)\bm w\big)
 +\alpha_\varepsilon(1-\phi)^2\bm w-\phi\nabla\pi
 &=\bm0 &&\text{in }D,\\
\nabla\cdot(\phi\bm w)&=0 &&\text{in }D,\\
\bm w=-\bm e_x\ \text{on }\Gamma_i\cup\Gamma_w,\qquad
\nu\partial_{\bm n}\bm w+(\bm u\cdot\bm n)\bm w+\pi\bm n
 &=-2u_x\bm e_x &&\text{on }\Gamma_o.
\end{aligned}
\end{equation}
Here $\bm e_x$ is the streamwise unit vector. The nonzero adjoint boundary
data account for the pressure and viscous terms in the drag objective;
inhomogeneous adjoint boundary data of this kind are the usual signature of a
force objective evaluated on the boundary \cite{giles2000introduction}.
For a design variation $\vartheta$ vanishing near $\partial D$, the sensitivity is
\begin{equation}
\label{eq:sens_flow}
\begin{aligned}
J'(\phi)[\vartheta]
= - \int_D \Big[&\nu\,\nabla\bm{u}:\nabla(\vartheta\,\bm{w})
 + \big((\bm{u}\cdot\nabla)(\vartheta\,\bm{u})\big)\cdot\bm{w}
 - 2\alpha_\varepsilon(1-\phi)\vartheta\,\bm{u}\cdot\bm{w}\\
 &- p\,\nabla\cdot(\vartheta\,\bm{w}) + \pi\,\nabla\cdot(\vartheta\,\bm{u}) \Big]\,\mathrm{d}\bm{x}.
\end{aligned}
\end{equation}
The state equation \eqref{eq:nse_ddm}, objective \eqref{eq:dragflow},
constraint \eqref{eq:volcon}, adjoint \eqref{eq:adjoint_flow}, and sensitivity \eqref{eq:sens_flow} form the
continuous optimality system. Appendix~\ref{app:flow} derives the adjoint and
sensitivity from the Lagrangian.

\paragraph{Discretization.}
We discretize \eqref{eq:nse_ddm} by a MAC staggered finite-difference scheme on a uniform grid over $D$, collecting velocities and pressures in $\mathbf U = (\bm{u}_h, p_h)$. Let $\bm{\phi}_f=\Pi_f\bm\phi$ denote the phase field averaged from cell centers to the velocity faces. The discrete state residual is
\begin{equation*}
\mathbf R(\mathbf U; \bm{\phi}) = K(\bm{\phi})\,\mathbf U + \mathbf N(\mathbf U;\bm{\phi}) - \mathbf F(\bm{\phi}) = 0,
\end{equation*}
with
\begin{equation}
\label{eq:Kphi_flow}
K(\bm{\phi}) = K_{\rm bc}
 + \mathrm{diag}(\bm{\phi}_f)\,K_{\rm mom}
 + \mathrm{diag}\!\Big(\alpha_\varepsilon(\bm{1}-\bm{\phi}_f)^2\Big)
 + K_{\rm div}\,\mathrm{diag}(\bm{\phi}_f),
\qquad \mathbf F(\bm{\phi}) = \mathbf F_{\rm bc} + \bm{\phi}_f \circ \mathbf F_{\rm mom},
\end{equation}
where $K_{\rm mom}$ carries the unweighted viscous and pressure-gradient stencils on the momentum rows, $K_{\rm div}$ the discrete divergence, $K_{\rm bc}$ the design-independent boundary and regularization terms, and $\mathbf N(\mathbf U;\bm\phi)$ is the discrete convection term $(\bm u\cdot\nabla)(\phi\bm u)$; the square of $\bm 1-\bm\phi_f$ is taken componentwise. The continuity rows include $-\delta_p p_h$, with a small constant $\delta_p>0$, to remove the pressure degeneracy in the solid; discrete weighted flux is therefore conserved up to this regularization. The penalty and momentum weighting act only on interior momentum rows. The four design-dependent contributions are the momentum row scaling, penalty diagonal, weighted divergence, and convection. Their coefficient structure permits exact reduced assembly; all four also enter the discrete sensitivity.


The state equations are solved by a damped Newton iteration with the exact Jacobian $J_N = K(\bm{\phi}) + \partial \mathbf N/\partial \mathbf U$, initialized with the state from the previous design. Let $J_h(\mathbf U)$ denote the discrete counterpart of \eqref{eq:dragflow}. It uses only values at $\partial D$, near which the design is fixed, so it depends on $\bm\phi$ only through $\mathbf U$. The discrete adjoint therefore solves
\begin{equation}
\label{eq:adjoint_disc_flow}
J_N^{\!\top}\, \mathbf W = \nabla_{\mathbf U} J_h ,
\end{equation}
and the gradient is $\nabla_{\bm{\phi}} J_h = -\,\Pi_f^{\!\top}(\partial \mathbf R/\partial\bm{\phi}_f)^{\!\top} \mathbf W$, where $\Pi_f^{\!\top}$ gathers face values back to the cells. Differentiating the design-dependent terms of $\mathbf R$ gives
\begin{equation*}
\Big(\frac{\partial \mathbf R}{\partial \bm{\phi}_f}\Big)^{\!\top} \mathbf W
= \mathbf W \circ \big( K_{\rm mom} \mathbf U - \mathbf F_{\rm mom} - 2\alpha_\varepsilon(\bm{1}-\bm{\phi}_f)\circ \mathbf U \big)\Big|_{\text{momentum rows}}
 \; + \; \big(K_{\rm div}^{\!\top} \mathbf W\big) \circ \mathbf U
 \; + \; \Big(\frac{\partial \mathbf N}{\partial \bm{\phi}_f}\Big)^{\!\top} \mathbf W .
\end{equation*}
The first term collects the momentum scaling, forcing and penalty, the second the weighted divergence, and the third the convection. Because \eqref{eq:adjoint_disc_flow} uses the transpose of the complete discrete Jacobian, including boundary rows and the pressure regularization, $\nabla_{\bm\phi}J_h$ is the gradient of the discrete objective under the implemented state equations.

\subsection{Compliance minimization of a cantilever}
\label{sec:elast}

The third problem minimizes the compliance of a cantilever, a structure governed by linear elasticity. It tests a symmetric state operator and a self-adjoint objective, which needs no separate adjoint solve. Let $\Omega \subset D$ be the region occupied by material, with $\partial\Omega = \Gamma_D \cup \Gamma_N \cup \Gamma_f$, where $\Gamma_D$ is clamped, $\Gamma_N$ carries a prescribed traction $\bm{t}$, and $\Gamma_f$ is the free design boundary (Figure~\ref{fig:settings}(c)). We denote the linearized strain by $\bm{S}(\bm{u}) = \tfrac12(\nabla\bm{u} + \nabla\bm{u}^{\!\top})$, since $\varepsilon$ is reserved for the interface width, and the linear elasticity operator that maps strain to stress by $A$, so that $\bm{\sigma}=A\,\bm{S}(\bm{u})$. For the isotropic material used here,
\begin{equation*}
A\,\bm{S}(\bm{u}) \;=\; 2\mu\,\bm{S}(\bm{u}) \;+\; \lambda\,\mathrm{tr}\big(\bm{S}(\bm{u})\big)\,\bm{I},
\end{equation*}
with $\lambda$ and $\mu$ the Lam\'e parameters. The displacement $\bm{u}$ satisfies
\begin{equation*}
\begin{gathered}
-\nabla\cdot\big(A\,\bm{S}(\bm{u})\big)=\bm{f}\qquad \text{in }\Omega,\\[-1mm]
\bm{u}=\bm{0}\ \text{on }\Gamma_D,\qquad
\big(A\,\bm{S}(\bm{u})\big)\bm{n}=\bm{t}\ \text{on }\Gamma_N,\qquad
\big(A\,\bm{S}(\bm{u})\big)\bm{n}=\bm{0}\ \text{on }\Gamma_f.
\end{gathered}
\end{equation*}

\paragraph{Diffuse-domain formulation.}
The evolving boundary is traction free, so we weight the stiffness in divergence form; no displacement penalty is needed there. The cantilever has fixed supports and loads; only its traction-free boundary moves with the design. We impose the clamp as an essential condition on $\Gamma_D$ and the load on the fixed $\Gamma_N$ (Figure~\ref{fig:settings}(c)). Following the smoothed-boundary mechanical model of \cite{yu2012extended} and the ersatz stiffness used in \cite{takezawa2010shape}, the weighted equation on the extended domain is
\begin{equation}
\label{eq:elast_ddm}
-\nabla\cdot\big(k(\phi)\,A\,\bm{S}(\bm{u})\big)
\;=\; \phi\,\bm{f}\qquad \text{in } D.
\end{equation}
Here $k(\phi)=k_{\min}+(1-k_{\min})\phi^p$ is the ersatz stiffness factor; it leaves the discrete operator affine in design-dependent coefficients. As in \cite{bendsoe2003topology,takezawa2010shape}, the power $p>1$ of the solid isotropic material with penalization (SIMP) model penalizes intermediate material fractions; the weight $k(\phi)$ of shape identification is linear in $\phi$. As the interface sharpens and $k_{\min}$ tends to zero, the weighted stress vanishes outside the material and the evolving boundary becomes traction free. With $V=\{\bm{v}\in[H^1(D)]^d:\bm{v}=\bm{0}\text{ on }\Gamma_D\}$, the weak problem actually solved is
\begin{equation}
\label{eq:elast_weak_design}
\int_D k(\phi)\, A\,\bm{S}(\bm{u}) : \bm{S}(\bm{v})\,\mathrm{d}\bm{x}
\;=\; \int_D \phi\,\bm{f}\cdot\bm{v}\,\mathrm{d}\bm{x} + \int_{\Gamma_N}\bm{t}\cdot\bm{v}\,\mathrm{d}s ,
\end{equation}
for every $\bm{v}\in V$. The positive floor makes the discrete operator nonsingular.

\paragraph{Optimality system.}
Following \cite{takezawa2010shape}, the experiments minimize the compliance together with a fixed multiple $f_V$ of the material volume, with no volume constraint,
\begin{equation}
\label{eq:compliance}
J(\phi) = J_{\rm c}(\phi) + f_V \int_D \phi\,\mathrm{d}\bm{x},
\qquad
J_{\rm c}(\phi) = \int_D \phi\,\bm{f}\cdot\bm{u} \,\mathrm{d}\bm{x} + \int_{\Gamma_N} \bm{t}\cdot\bm{u}\,\mathrm{d}s .
\end{equation}
The compliance is self-adjoint. In weak form its adjoint $\bm{w}$ solves
\begin{equation*}
\int_D k(\phi)A\bm{S}(\bm{w}):\bm{S}(\bm{v})\,\mathrm{d}\bm{x}
=\int_D \phi\bm{f}\cdot\bm{v}\,\mathrm{d}\bm{x}
 +\int_{\Gamma_N}\bm{t}\cdot\bm{v}\,\mathrm{d}s
\quad\text{for every }\bm{v}\in V,
\end{equation*}
which is the state equation \eqref{eq:elast_weak_design}; hence
$\bm{w}=\bm{u}$ and no additional solve is required. For a load supported
where the design is fixed, the sensitivity
$J'(\phi)[\vartheta]=\int_D J'(\phi)\,\vartheta\,\mathrm{d}\bm{x}$ has the
$L^2(D)$ representative
\begin{equation}
\label{eq:sens_compliance}
J_{\rm c}'(\phi) = -\,k'(\phi)\,A\,\bm{S}(\bm{u}) : \bm{S}(\bm{u}),
\qquad
J'(\phi) = J_{\rm c}'(\phi) + f_V .
\end{equation}
Appendix~\ref{app:elast} derives the adjoint and sensitivity from the
Lagrangian.

\paragraph{Discretization.}
We discretize displacement with bilinear ($Q_1$) finite elements on a uniform quadrilateral grid of mesh size $h$ and take one piecewise-constant design value $\phi_e$ per element. The stiffness matrix is
\begin{equation}
\label{eq:elast_affine}
K(\bm{\phi}) \;=\; \sum_{e} k(\phi_e)\,K_e ,
\end{equation}
with design-independent element matrices $K_e$; the state solves $K(\bm{\phi})\,\mathbf{U} = \mathbf{F}$. The discrete objective is $J_h=\mathbf{F}^{\!\top}\mathbf{U}+f_Vh^2\,\bm 1^{\!\top}\bm\phi$. Its compliance part is self-adjoint, so the discrete adjoint equals the state, $\mathbf{W} = \mathbf{U}$, and the gradient has the entries
\[
(\nabla_{\bm\phi}J_h)_e = -\,k'(\phi_e)\,\mathbf{U}_e^{\!\top} K_e\, \mathbf{U}_e + f_V h^2 .
\]
Here $\mathbf{U}_e$ contains the displacements of element $e$.

\subsection{The gated Allen--Cahn design update}
\label{sec:design}

Once the state and adjoint are available, their sensitivity drives the
interface toward a better design, while diffusion and the double-well term
maintain a transition layer. We use the gated Allen--Cahn update developed for the companion
Brinkman study \cite{li2026residual}; only the ingredients needed to
reproduce the three DDM runs are stated here. With the double-well
potential $F(\phi)=\tfrac14\phi^2(1-\phi)^2$, the gate $G'(\phi)=\phi^2(1-\phi)^2$,
and $\hat\varepsilon=\varepsilon/(6\sqrt2)$, the design flow is
\begin{equation}
\label{eq:acflow}
\partial_t\phi=\hat\varepsilon^2\Delta\phi-F'(\phi)
-\eta G'(\phi)\widehat{J'} .
\end{equation}
Here $J'(\phi)$ is the $L^2(D)$ representative of the sensitivity; the
computations use the discrete gradient $\nabla_{\bm\phi}J_h$ in its place.
The choice of $\hat\varepsilon$ makes the stationary profile of
\eqref{eq:acflow} coincide with \eqref{eq:phase}; the different
interface-width convention used in \cite{li2026residual} accounts for
the apparent factor of two in its coefficient. The gate confines design
forcing to the interface. All three problems use the same after-gate normalization:
\begin{equation}
\label{eq:forcing}
\widehat{J'}=\frac{J'(\phi)}{\|G'(\phi)J'(\phi)\|_{\rm disc}+\delta},
\qquad q(\phi)=G'(\phi)\widehat{J'},
\end{equation}
where $\|\cdot\|_{\rm disc}$ is the unweighted Euclidean norm of the stored
sensitivity array and a small $\delta>0$ prevents division by zero.
To choose the forcing strength by the same rule in every problem, we
compute the initial full-order sensitivity and set
\begin{equation}
\label{eq:balance}
\eta=\rho\,\frac{\|F'(\phi_0)\|_{B_0}}{\|q(\phi_0)\|_{B_0}},
\qquad B_0=\{0.1<\phi_0<0.9\}.
\end{equation}
Here $\rho>0$ is a balance factor, the norms are restricted to the initial
transition band $B_0$, and $\eta$ is then held fixed. This calibrates the initial forcing relative to the
double-well term, not the complete free-energy derivative, whose diffusion
and potential terms can nearly cancel. In the volume-constrained case the
ratio is measured before the multiplier correction. We note that it is a common
amplitude-selection rule, not an energy-decay or convergence guarantee.
The gate cannot nucleate a hole in a region with exactly binary $\phi$.

For the flow problem, the constraint \eqref{eq:volcon} is enforced by adding
$\lambda G'(\phi)$ to \eqref{eq:acflow}:
\begin{equation*}
\partial_t\phi=\hat\varepsilon^2\Delta\phi-F'(\phi)
-\eta G'(\phi)\widehat{J'}+\lambda G'(\phi).
\end{equation*}
The cantilever instead includes a volume penalty in
\eqref{eq:compliance}; shape identification has no volume condition.
At each design iteration we take $N_{\rm in}$ semi-implicit steps of size $\Delta t$.
The sensitivity stays fixed throughout these inner steps. Writing
$M=I-\Delta t\,\hat\varepsilon^2L_h$, with $L_h$ the discrete Laplacian, and
$\bm r_0=\bm\phi+\Delta t\,[-F'(\bm\phi)-
\eta G'(\bm\phi)\circ\widehat{\nabla_{\bm\phi}J_h}]$, the unconstrained
and volume-constrained steps are respectively
\begin{equation}
\label{eq:ac_update}
M\bm\phi^{\,\rm new}=\bm r_0,
\end{equation}
\begin{equation}
\label{eq:ac_update_vol}
M\bm\phi^{\,\rm new}=\bm r_0+\lambda\bm r_1,
\qquad \bm r_1=\Delta t\,G'(\bm\phi).
\end{equation}
Let $M\bm\phi_a=\bm r_0$ and $M\bm\phi_b=\bm r_1$. The multiplier is
chosen at each inner step to enforce the discrete volume exactly:
\begin{equation}
\label{eq:lambda}
\lambda=\frac{V_0/h^d-\bm 1^{\!\top}\bm\phi_a}
{\bm 1^{\!\top}\bm\phi_b}.
\end{equation}
Here projection onto $[0,1]$ means componentwise clipping of the updated design
to $[0,1]$. When the clipping is active, that is, when $\bm\phi_a+\lambda\bm\phi_b$
leaves $[0,1]$, the quotient \eqref{eq:lambda} no longer gives the exact volume, and
$\lambda$ is instead found by bisection so that the clipped field
$\min(\max(\bm\phi_a+\lambda\bm\phi_b,0),1)$ has the prescribed volume.
An initial design must have an interfacial band; when $\phi$ is binary,
$G'(\phi)$ and the denominator vanish. We initialize with the signed-distance profile \eqref{eq:phase}. The design boundary conditions,
volume treatments, projections, and solve domains used in each problem are
listed in Table~\ref{tab:design}. These are the only design-update
differences relevant to the ROM.

\begin{table}[htbp]
\centering
\caption{Design-update choices used in the reported runs for the three DDM problems. For shape identification the semi-implicit solve is restricted to the interfacial band; drag minimization and compliance minimization use the full grid.}
\label{tab:design}
\small\setlength{\tabcolsep}{4pt}
\begin{tabular}{llllll}
\toprule
problem & design values & $\phi$ on $\partial D$ & volume & projection & solve domain \\
\midrule
shape identification & nodes & $\phi=0$ & free & no & interface band \\
drag minimization & cells & $\phi=1$ & multiplier \eqref{eq:lambda} & yes & full grid \\
compliance minimization & elements & Neumann & penalty in \eqref{eq:compliance} & yes & full grid \\
\bottomrule
\end{tabular}
\end{table}

The normalized forcing can leave a nonzero design-step floor even after
the objective has settled, while an oscillatory design can have a flat
windowed objective. We therefore use the objective test
\begin{equation}
\label{eq:stoptest}
|J_k-J_{k-N_{\rm stop}}|
\le \mathrm{tol}_J\,N_{\rm stop}\,|J_k|
\end{equation}
at two successive iterations, with window length $N_{\rm stop}$ and tolerance
$\mathrm{tol}_J$; here $J_k$ denotes $J_h$ at the $k$th design, evaluated
with the state used for that design update. Independently, we stop if the maximum
$\|\bm\phi_{j+1}-\bm\phi_j\|_\infty$ over the last $N_{\rm stop}$ updates is
below $\mathrm{tol}_\phi$. Either test describes stagnation, not stationarity;
Appendix~\ref{sec:tau} examines sensitivity to the design step. This completes the full-order optimization: at each design iteration, the state and adjoint equations of Sections~\ref{sec:opt}--\ref{sec:elast} are solved at full order, the sensitivity is evaluated, and the design is advanced by \eqref{eq:ac_update} or \eqref{eq:ac_update_vol}.

Full-order state and adjoint solves dominate the cost of each iteration. Next, we present CL-ROM, which replaces them by local reduced-order solves whenever these pass its residual and objective checks.

\FloatBarrier
\section{The Compact Residual-Screened Local Reduced-Order Model}
\label{sec:trrom}

We retain the decision loop of CL-ROM \cite{li2026residual}, summarized in
Algorithm~\ref{alg:trrom}. Recent full-order and corrected solutions span
small POD spaces. At each new design, the reduced state and adjoint computed
in these spaces are accepted if the residuals of the full discrete equations
pass a screen; otherwise they are corrected with retained factors or replaced
by a full-order solve, and the objective is verified at full order when its estimated error is large.
Section~\ref{sec:algorithm} states the ingredients of the algorithm briefly,
referring to \cite{li2026residual} for their motivation, and gives a local
approximation estimate for the DDM systems. Section~\ref{sec:models}
describes the operators and solvers that depend on the problem.

\subsection{The algorithm and its accuracy}
\label{sec:algorithm}

\begin{algorithm}[htbp]
\small
\caption{CL-ROM for all three diffuse-domain design problems, adapted from \cite{li2026residual}}
\label{alg:trrom}
\begin{algorithmic}[1]
\State \textbf{Input:} initial design $\bm{\phi}_0$; screening tolerance $\tau_\eta \in [\tau_\eta^{\min}, \tau_\eta^{\max}]$, objective tolerance $\tau_J$, thresholds $\theta_r < \theta_v < 1$; window and rank caps $w_{\max}$, $r_{\max}$; warm-up $N_{\rm w}$, relaxation period $N_c$, backstop $N_{\max}$, refactoring threshold $m_{\max}$
\State $N_{\rm w}$ full-order design iterations; retain the latest factors of $K_{\rm ref}$; fill the state and adjoint windows; POD bases $\Psi_s$, $\Psi_a$
\For{$k = N_{\rm w}, N_{\rm w}+1, \dots$ until the stopping rule of Section~\ref{sec:design} triggers}
    \State reduced state and adjoint in the spans of $\Psi_s$ and $\Psi_a$ (Section~\ref{sec:models}); calibrated indicators $\tilde\eta_s, \tilde\eta_a$ by \eqref{eq:etacal}; objective estimate $\eta_J$ by \eqref{eq:dwr}
    \If{$\max(\tilde\eta_s, \tilde\eta_a) > \tau_\eta$}
        \State try up to two corrections \eqref{eq:compact_correction} (one for flow), screening each corrected pair
    \EndIf
    \If{the candidate passes the residual screen}
        \If{$|\eta_J| < \theta_r\tau_J|J_h|$ throughout the last $N_c$ accepted iterations}
            \State $\tau_\eta \gets \min(2\tau_\eta, \tau_\eta^{\max})$
        \EndIf
        \If{$|\eta_J| > \theta_v\tau_J|J_h|$ and the backstop has not been reached}
            \State try up to two further corrections; require both the residual and objective checks to pass
        \EndIf
    \EndIf
    \If{screening fails, the objective check still fails, or the backstop is reached}
        \State full-order pair with the solver of Table~\ref{tab:rommodels}; replace the candidate and update calibration
        \State at verification, adjust $\tau_\eta$ using the observed objective error
    \EndIf
    \State sensitivity at the pair obtained; design update \eqref{eq:ac_update} or \eqref{eq:ac_update_vol}
    \State append any full-order or accepted corrected pair to the windows and rebuild the bases
\EndFor
\end{algorithmic}
\end{algorithm}

\paragraph{Windows and bases.}
Separate state and adjoint windows retain at most $w_{\max}$ recent
full-order or accepted corrected snapshots. An uncorrected reduced solution
already lies in its current space and does not enrich that space.
POD gives bases $\Psi_s$ and $\Psi_a$ with at most
$r_{\max}$ modes above a relative singular-value cutoff (the compliance objective \eqref{eq:compliance} needs only the state basis). Depending on the model, the reduced systems use
Galerkin or least-squares Petrov--Galerkin projection
(Section~\ref{sec:models}). Enrichment updates the bases and
basis-dependent operator products, while sensitivities use the
problem's full-order formula evaluated at the reduced fields.

The a priori estimate in \cite{li2026residual} separates POD truncation from
the effect of design changes, which enter that model only through a diagonal
friction term. For the DDM systems, moving from a snapshot design
$\bm\phi_j$ to the current design changes $A(\bm\phi)$ through $k(\phi)$; for shape
identification it also changes $A(\bm\phi)$ through $\delta_\sigma(\phi)$ and
the load $\mathbf b(\bm\phi)$ through $\phi f$. Corrected snapshots
also have a nonzero residual. The next result makes these contributions
explicit for the two symmetric models, with all essential conditions
eliminated and $\|\mathbf v\|_B=(\mathbf v^\top B\mathbf v)^{1/2}$ for a
symmetric positive definite (SPD) matrix $B$.

\begin{samepage}
\begin{proposition}[Local approximation with inexact snapshots]
\label{prop:ddm_local}
Let $A\mathbf U=\mathbf b$ be the current linear, symmetric positive definite DDM system (as in the shape-identification and cantilever problems) and let $\mathbf U_r$
be its Galerkin solution in the space $V_r$ spanned by the state basis. At snapshot design $\bm\phi_j$, write
$A_j=A(\bm\phi_j)$, $\mathbf b_j=\mathbf b(\bm\phi_j)$, and let $\widetilde{\mathbf U}_j$
be the possibly corrected field used to form the snapshot. Suppose
$m_j,M_j>0$ and $m_j A_j\preceq A\preceq M_j A_j$, where $\preceq$ denotes the
ordering of quadratic forms. Set
\begin{equation*}
\mathbf s_j=\mathbf b_j-A_j\widetilde{\mathbf U}_j,\qquad
\mathbf d_j=(\mathbf b-\mathbf b_j)-(A-A_j)\widetilde{\mathbf U}_j,\qquad
t_j=\|(I-P_r)\widetilde{\mathbf U}_j\|_{A_j},
\end{equation*}
where $P_r$ is the Euclidean orthogonal projection onto $V_r$. Then
\begin{equation}
\label{eq:ddm_local_bound}
\|\mathbf U-\mathbf U_r\|_A\le
\min_j\left\{
\frac{\|\mathbf d_j\|_{A_j^{-1}}+\|\mathbf s_j\|_{A_j^{-1}}}{\sqrt{m_j}}
+\sqrt{M_j}\,t_j\right\}.
\end{equation}
\end{proposition}
\end{samepage}
\begin{proof}
Galerkin orthogonality gives $\|\mathbf U-\mathbf U_r\|_A=\inf_{\mathbf v\in V_r}\|\mathbf U-\mathbf v\|_A$.
For any stored snapshot, $A(\mathbf U-\widetilde{\mathbf U}_j)=\mathbf s_j+\mathbf d_j$, so
$\|\mathbf U-\widetilde{\mathbf U}_j\|_A=\|\mathbf s_j+\mathbf d_j\|_{A^{-1}}
\le m_j^{-1/2}\big(\|\mathbf s_j\|_{A_j^{-1}}+\|\mathbf d_j\|_{A_j^{-1}}\big)$,
since $A\succeq m_jA_j$ implies $A^{-1}\preceq m_j^{-1}A_j^{-1}$.
Taking $\mathbf v=P_r\widetilde{\mathbf U}_j$ and
$\|(I-P_r)\widetilde{\mathbf U}_j\|_A\le\sqrt{M_j}\,t_j$ proves the claim.
\end{proof}

The three terms of \eqref{eq:ddm_local_bound} measure operator and load changes, snapshot error, and POD
truncation. Recent snapshots reduce $\mathbf d_j$ and accurate corrections
reduce $\mathbf s_j$, but neither a short window nor discarding old
snapshots guarantees a small truncation term, which the residual screen must
detect. The models differ in the first term: for shape identification,
$A-A_j$ includes the change of $\delta_\sigma(\phi)$ and hence of discrete
phase-field gradients, whereas for elasticity it depends on the design only
through elementwise stiffness ratios, and the load is fixed.

\paragraph{Screening and its calibration.}
\label{sec:indicators}
The residuals of a reduced pair $(\mathbf U_r,\mathbf W_r)$ are $\mathbf r_s=\mathbf b-A\mathbf U_r$ and $\mathbf r_a=\nabla_{\mathbf U}J_h(\mathbf U_r)-A^{\!\top}\mathbf W_r$ for the linear problems, and $\mathbf r_s=-\mathbf R(\mathbf U_r;\bm\phi)$ and $\mathbf r_a=\nabla_{\mathbf U}J_h(\mathbf U_r)-J_N(\mathbf U_r)^{\!\top}\mathbf W_r$ for the flow, all evaluated in the full discrete space. Their relative norms $\eta^{E}_s=\|\mathbf r_s\|/\|\mathbf b(\bm\phi)\|$ and $\eta^{E}_a=\|\mathbf r_a\|/\|\nabla_{\mathbf U}J_h(\mathbf U_r)\|$ measure error in the discrete equations at a fixed design, but their relation to the solution error depends on the unknown conditioning of the operator. As in \cite{li2026residual}, they are therefore calibrated against the preconditioned indicators
\begin{equation}
\label{eq:etaP}
\eta_s = \frac{\|K_{\rm ref}^{-1}\mathbf{r}_s\|}{\|\mathbf{U}_{\rm ref}\|},
\qquad
\eta_a = \frac{\|K_{\rm ref}^{-\top}\mathbf{r}_a\|}{\|K_{\rm ref}^{-\top}\nabla_{\mathbf U}J_h(\mathbf{U}_r)\|},
\end{equation}
formed using the factors of $K_{\rm ref}$ retained from the most recent refactorization and the state $\mathbf{U}_{\rm ref}$ computed at that design. Since applying the factors at every candidate would add full-dimensional triangular solves, acceptance uses the \emph{calibrated indicators}
\begin{equation}
\label{eq:etacal}
\tilde\eta_s = \eta^{E}_s / c_s, \qquad \tilde\eta_a = \eta^{E}_a / c_a,
\end{equation}
where $c_s$ and $c_a$ are the medians of $\eta^{E}_s/\eta_s$ and $\eta^{E}_a/\eta_a$ over the ten most recent evaluations of \eqref{eq:etaP}, made at the first three reduced iterations and at every fallback and verification; corrected pairs keep separate histories, used once three ratios are available. A pair passes when $\max(\tilde\eta_s, \tilde\eta_a) \le \tau_\eta$. The threshold applies to a calibrated residual and is not an error tolerance.

\paragraph{Objective-error indicator.}
A small state residual need not imply a small objective error. As in \cite{li2026residual}, the objective error is estimated by a dual-weighted residual with an enriched dual,
\begin{equation}
\label{eq:dwr}
\eta_J=\mathbf W_+^{\!\top}\mathbf r_s,
\qquad \mathbf W_+\in\mathrm{span}\{\Psi_a,\mathbf W_r,\mathbf r_s\},
\end{equation}
where the residual direction is orthogonalized against the adjoint space, $\mathbf W_r$ adds a direction only when correction has moved it outside that space, and $\mathbf W_+$ is obtained by least-squares projection for the flow and by Galerkin projection for the symmetric problems. It approximates the leading term $\mathbf W^{\!\top}\mathbf r_s$ of the objective error, with $\mathbf W$ the full-order adjoint, up to a remainder quadratic in the state error \cite{li2026residual}. The estimate is inexpensive to evaluate, although it is not an upper bound.

\paragraph{Corrections, fallbacks and verification.}
A rejected linear pair is corrected with the retained factors,
\begin{equation}
\label{eq:compact_correction}
\mathbf{U}^{+} = \mathbf{U}_r + K_{\rm ref}^{-1}\,\mathbf{r}_s(\mathbf{U}_r), \qquad
\mathbf{W}^{+} = \mathbf{W}_r + K_{\rm ref}^{-\top}\,\mathbf{r}_a(\mathbf{W}_r; \mathbf{U}^{+}).
\end{equation}
For the state equation, this is one preconditioned Richardson step; it
reduces the error when the preconditioned error propagator
$I-K_{\rm ref}^{-1}A$ is contractive in the relevant norm. The algorithm
does not test this condition but screens the corrected pair again by its
full residual. For the
nonlinear flow, a correction is one inexact Newton step. If screening still
fails, a full-order solve starts from the candidate with the solver of
Table~\ref{tab:rommodels}. At verification, $\tau_\eta$ is halved if the observed relative
objective error exceeds $\tau_J$ and doubled if it is below
$\theta_r\tau_J$, within the prescribed bounds. Full-order and
accepted corrected pairs both enter the snapshot windows, so a correction
adds information outside the current reduced spaces even when no full-order
solve is made.

\subsection{What differs between the problems}
\label{sec:models}

Table~\ref{tab:rommodels} lists the ingredients of Algorithm~\ref{alg:trrom} that depend on the problem. The solver choices follow two distinctions. Symmetric positive definite systems admit Galerkin projection and Cholesky-based solves; the nonlinear saddle-point flow system uses least-squares projection and Newton--GMRES. A self-adjoint compliance objective needs only the state, whereas the other implementations retain separate state and adjoint spaces. In all three problems the design enters through pointwise functions of $\bm{\phi}$ multiplying fixed sparse stencils, so the reduced operators are assembled exactly from the current coefficients and basis vectors without empirical interpolation; evaluating the exact full-order residual of a reduced pair still costs $O(nr)$ sparse and elementwise work, where $n$ is the number of full-order unknowns and $r$ the reduced dimension.

\begin{table}[htbp]
\centering
\caption{Problem-dependent ingredients of CL-ROM. $K_{\rm ref}$ is the matrix whose factors are retained. The fallback solver refactors when its Krylov iteration count exceeds $m_{\max}$; for drag minimization the count runs over the Newton and adjoint solves of an update. After an objective trigger, all three problems allow up to two corrections.}
\label{tab:rommodels}
\small\setlength{\tabcolsep}{4pt}
\begin{tabular}{@{}>{\centering\arraybackslash}p{0.19\textwidth}>{\raggedright\arraybackslash}p{0.18\textwidth}>{\raggedright\arraybackslash}p{0.285\textwidth}>{\raggedright\arraybackslash}p{0.285\textwidth}@{}}
\toprule
 & shape identification & drag minimization & compliance minimization \\
\midrule
state operator & SPD, \eqref{eq:state_discrete} & nonlinear saddle point, \eqref{eq:Kphi_flow} & SPD, \eqref{eq:elast_affine} \\
objective & not self-adjoint & not self-adjoint & self-adjoint \\
reduced problems & \eqref{eq:rom_scalar}, \eqref{eq:rom_assembly} & \eqref{eq:rom_flow}, \eqref{eq:rom_apply_flow} & \eqref{eq:rom_elast} \\
reduced solves & Galerkin, two bases & LSPG with Gauss--Newton, two bases & Galerkin, one basis \\
indicators screened & state, adjoint & state, adjoint & state \\
retained factors & Cholesky of $A(\bm{\phi})$ & LU of the Newton Jacobian & Cholesky of $K(\bm{\phi})$ \\
fallback & CG & Newton--GMRES & CG \\
correction after a failed screen & \eqref{eq:compact_correction}, two attempts & inexact Newton step, one attempt & \eqref{eq:compact_correction} on the state, two attempts \\
\bottomrule
\end{tabular}
\end{table}

\paragraph{Shape identification.}
The operator \eqref{eq:state_discrete} is symmetric positive definite. The reduced state and adjoint use Galerkin projection onto separate bases, since the objective gives the adjoint a different load:
\begin{equation}
\label{eq:rom_scalar}
\begin{aligned}
\Psi_s^{\!\top}A(\bm\phi)\Psi_s\,\mathbf a&=\Psi_s^{\!\top}\mathbf b(\bm\phi), &\qquad \mathbf U_r&=\Psi_s\mathbf a,\\
\Psi_a^{\!\top}A(\bm\phi)\Psi_a\,\mathbf c&=\Psi_a^{\!\top}\nabla_{\mathbf U}J_h(\mathbf U_r), &\qquad \mathbf W_r&=\Psi_a\mathbf c .
\end{aligned}
\end{equation}
Precomputing $G\Psi_s$ and $G\Psi_a$ at each basis update permits the operator actions
\begin{equation}
\label{eq:rom_assembly}
A(\bm\phi)\Psi=G^\top\!\left(\big[k_{\min}+(1-k_{\min})\Pi_f\bm\phi\big]\circ(G\Psi)\right)
+h^2\lambda_\varepsilon\,\delta_\sigma(\bm\phi)\circ\Psi,
\end{equation}
where $\Psi$ is either basis and vector factors act columnwise. These actions use sparse edge operations, without constructing a separate full matrix for each basis vector; the residuals of the reduced pair are still assembled in the full discrete space for screening. Each basis action takes $O(nr)$ sparse arithmetic, forming a projected matrix takes $O(nr^2)$, and each small solve costs $O(r^3)$. The reduced sensitivity is \eqref{eq:sensitivity} at the reduced state and adjoint.

\paragraph{Drag minimization in a channel.}
The flow system is nonlinear and nonsymmetric, so the reduced state and adjoint minimize the full-order residuals over their bases (least-squares Petrov--Galerkin, LSPG \cite{carlberg2011efficient}):
\begin{equation}
\label{eq:rom_flow}
\begin{aligned}
\mathbf a&=\operatorname*{arg\,min}_{\mathbf a}\,\big\|\mathbf R(\Psi_s\mathbf a;\bm\phi)\big\|_2, &\qquad \mathbf U_r&=\Psi_s\mathbf a,\\
\mathbf c&=\operatorname*{arg\,min}_{\mathbf c}\,\big\|J_N(\mathbf U_r)^{\!\top}\Psi_a\mathbf c-\nabla_{\mathbf U}J_h(\mathbf U_r)\big\|_2, &\qquad \mathbf W_r&=\Psi_a\mathbf c .
\end{aligned}
\end{equation}
The state problem is solved by Gauss--Newton, starting from the projection of the last accepted state, with increments $\delta\mathbf a$ that minimize $\|J_N(\Psi_s\mathbf a)\Psi_s\,\delta\mathbf a+\mathbf R(\Psi_s\mathbf a;\bm\phi)\|_2$. The design enters three linear blocks of $K(\bm\phi)$, each as an elementwise factor multiplying a fixed operator, and the convection term. For a basis $\Psi$,
\begin{equation}
\label{eq:rom_apply_flow}
\begin{aligned}
K(\bm{\phi})\,\Psi &= K_{\rm bc}\Psi + \bm{\phi}_f \circ (K_{\rm mom}\Psi)
 + \alpha_\varepsilon(\bm{1}-\bm{\phi}_f)^2 \circ \Psi + K_{\rm div}\big(\bm{\phi}_f \circ \Psi\big),\\
K(\bm{\phi})^{\!\top}\Psi &= K_{\rm bc}^{\!\top}\Psi + K_{\rm mom}^{\!\top}\big(\bm{\phi}_f \circ \Psi\big)
 + \alpha_\varepsilon(\bm{1}-\bm{\phi}_f)^2 \circ \Psi + \bm{\phi}_f \circ \big(K_{\rm div}^{\!\top}\Psi\big).
\end{aligned}
\end{equation}
The fixed operator--basis products are reused after each basis update, while the convection and its Jacobian are re-evaluated at each Gauss--Newton step. A pair rejected by the screen receives one inexact Newton correction, with Newton--GMRES as the fallback.

\paragraph{Compliance minimization of a cantilever.}
Compliance is self-adjoint, so a single Galerkin basis $\Psi$ serves state and adjoint. Projecting \eqref{eq:elast_affine} gives
\begin{equation}
\label{eq:rom_elast}
K_r(\bm\phi)\,\mathbf a=\Psi^{\!\top}\mathbf F,\qquad
K_r(\bm\phi)=\sum_e k(\phi_e)\,\Psi^{\!\top}K_e\Psi,\qquad
\mathbf U_r=\mathbf W_r=\Psi\mathbf a .
\end{equation}
Each $r\times r$ product $\Psi^{\!\top}K_e\Psi$ involves only the rows of $\Psi$ on element $e$ and is precomputed at each basis update, so a new design requires only their weighted sum. Only the state indicator is screened. For an uncorrected Galerkin state, Galerkin orthogonality makes the compliance error equal to $\mathbf r_s^{\!\top}K^{-1}\mathbf r_s$, which is quadratic in the state residual. A rejected cantilever pair receives up to two corrections \eqref{eq:compact_correction} of the state, with CG preconditioned by the retained Cholesky factors as the fallback, under the same refresh rule as the full-order baseline.

\FloatBarrier
\section{Numerical Experiments}
\label{sec:experiments}

The experiments address two questions: does the reduced solver preserve the objective and design trajectories of the full-order solver, and how much work does it save? For each problem, the full-order and reduced runs use the same DDM construction and design update, changing only the state--adjoint solver. The comparisons therefore test transfer of CL-ROM across models, rather than selecting a best diffuse formulation. The three problem studies make these comparisons; supporting tests are collected in the appendices.

\subsection{Common settings}
\label{sec:protocol}

CL-ROM is tested on the three formulations of
Section~\ref{sec:formulation}: shape identification, drag minimization in a channel,
and compliance minimization of a cantilever. We report one grid for each main ROM comparison;
Appendix~\ref{app:ddm_refinement} checks mesh refinement of each fixed-geometry
diffuse model at the corresponding interface width.
All three use $\varepsilon=4h$,
$\Delta t=0.01$, $N_{\rm in}=5$, the normalization \eqref{eq:forcing} with
$\delta=10^{-14}$, and the initial balance rule \eqref{eq:balance} with
$\rho=0.3$. The resulting fixed coefficients
are $\eta=0.70$, $0.13$, and $0.72$, respectively. No problem-specific forcing sweep is used to select these coefficients, and within each problem the full-order and CL-ROM runs use the same initial design and coefficient. The stopping test
\eqref{eq:stoptest} uses $N_{\rm stop}=50$, $\mathrm{tol}_J=10^{-6}$ and
$\mathrm{tol}_\phi=10^{-5}$.

Unless stated otherwise, CL-ROM starts with $N_{\rm w}=2$ full-order updates,
uses window and rank caps $w_{\max}=r_{\max}=8$, and retains POD
modes above a relative cutoff of $10^{-10}$. For the nonlinear flow, the
reduced state solve uses at most two Gauss--Newton steps per design iteration.
The screening tolerance starts
at $\tau_\eta=10^{-5}$ and remains in
$[\tau_\eta^{\min},\tau_\eta^{\max}]=[10^{-9},3\times10^{-5}]$.
Objective verification uses $\tau_J=10^{-4}$, $\theta_v=1/2$,
$\theta_r=1/10$, $N_c=25$, and $N_{\max}=400$. Full-order Krylov solves use
relative residual tolerance $10^{-10}$. The symmetric problems refresh
their retained factors beyond $m_{\max}=8$ iterations; the nonlinear flow
counts iterations over its Newton and adjoint solves and uses $m_{\max}=40$.

A larger design step, $\Delta t=0.05$, was also tested on all three problems.
It reduced the optimization
time for the flow and the cantilever while retaining comparable
objective levels and resolved interfaces on those grids. For shape identification,
however, the larger step overshot the lowest objective, increased the
geometric error, and gave little reduction in CL-ROM time. We therefore
retain the conservative common value $\Delta t=0.01$ in the main comparisons.
Appendix~\ref{sec:tau} reports the sensitivity results for all three cases.

As in the companion Brinkman study \cite{li2026residual}, the full-order
baseline, labeled FOM in the tables, solves the state and adjoint at every design with the Krylov
solvers above, preconditioned by retained factors that are refreshed only
when the iteration count exceeds $m_{\max}$. CL-ROM uses the same
full-order solver and refresh rule for fallbacks and verifications. The
baseline's stopping test determines the comparison length, and CL-ROM is
timed for this same number of updates. Separate runs re-solve
residual-screened candidates at full order at fixed designs
(Section~\ref{sec:audit}); their timings are excluded from
speedup comparisons.

The main runs were timed sequentially in single-threaded MATLAB on one
Linux workstation with an Intel Core i9-10900F processor.\footnote{The authors used OpenAI Codex and Anthropic Claude Code to assist with code development, reviewing and editing the generated outputs as needed. The code is available at \url{https://github.com/zwang01/diffuse_domain_CLROMs}.} Wall times include
setup and the initial coefficient calibration, but exclude post-run
full-order scoring. Solve and factorization counts refer to the optimization
loop. These are individual timing measurements. Tables report mean objectives
over the last fifty updates, and relative differences from the corresponding
full-order mean.

\subsection{Shape identification}
\label{sec:num}

Shape identification tests whether reduced solves preserve recovery of a prescribed geometry. We test Algorithm~\ref{alg:trrom} on the shape-identification problem \eqref{eq:objective} on $D = (-0.15, 1.25) \times (-0.1, 0.1)$. The reference domain $\Omega_{\rm ref}$, of area $0.0693$, is enclosed by $\big(x(t),y(t)\big) = \big((t-1)^2,\,0.13t(t-1)(t-2)\big)$, $t\in[0,2]$. Its airfoil-like profile has a cusp at the trailing edge. We use $f=100(0.1014x-2.0676)$, $g_D=0$, and the sharp-domain target $\bar u=100(y^2-0.0169x(x-1)^2)\chi_{\Omega_{\rm ref}}$. This example is taken from \cite{laporte2003numerical}. The surface penalty uses $c_D=0.1$ and $\sigma=10^{-5}/\varepsilon$, and the diffusion-coefficient floor is $k_{\min}=10^{-6}$. The initial design is a thin heptagon. The reported CL-ROM comparison uses $h=2\times10^{-3}$ ($n=69{,}201$).

For this problem, the semi-implicit design solve is restricted to
$\{\phi(1-\phi)>10^{-12}\}$ enlarged by two cells; all compared runs use
this same update. Full-order solves use CG preconditioned by retained
Cholesky factors.

\begin{figure}[htbp]
\centering
\includegraphics[width=0.49\textwidth]{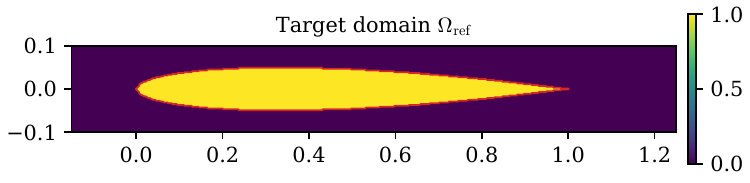}\hfill
\includegraphics[width=0.49\textwidth]{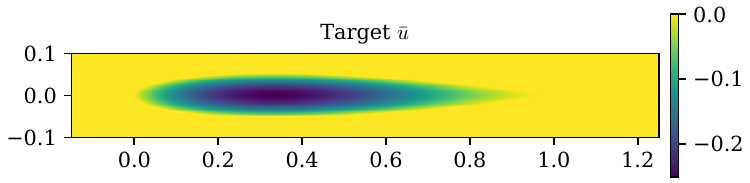}\\[1mm]
\includegraphics[width=0.49\textwidth]{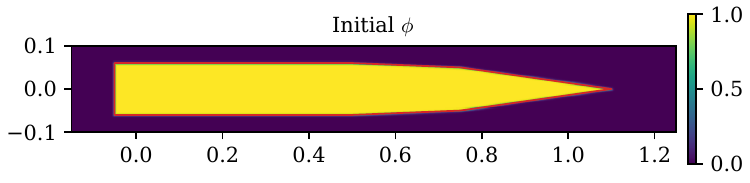}\hfill
\includegraphics[width=0.49\textwidth]{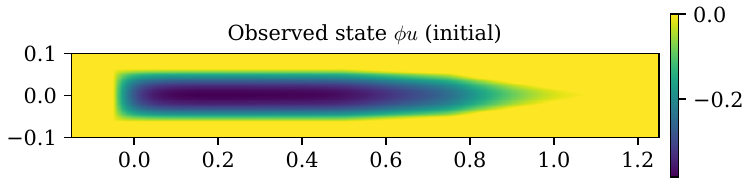}\\[1mm]
\includegraphics[width=0.49\textwidth]{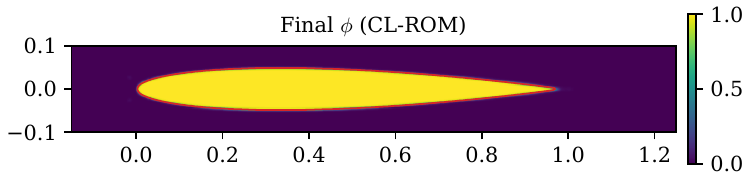}\hfill
\includegraphics[width=0.49\textwidth]{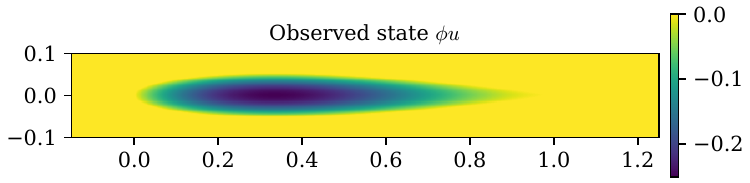}
\caption{Shape identification at $h=2\times10^{-3}$, panels at true aspect ratio, with the $\phi=\tfrac12$ contour in red. Top: the target domain $\Omega_{\rm ref}$ and the target state $\bar u$. Middle: the initial phase field and its observed state $\phi u$, whose color scale extends below that of the other two states. Bottom: the CL-ROM design after 7139 updates and its observed state.}
\label{fig:phase}
\end{figure}
\begin{figure}[htbp]
\centering
\includegraphics[width=0.49\textwidth]{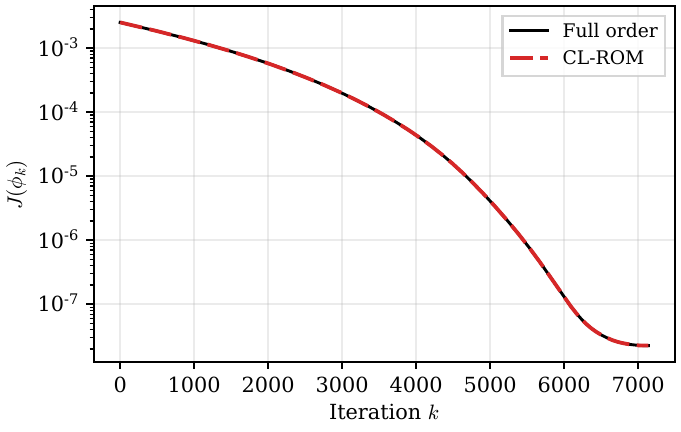}\hfill
\includegraphics[width=0.49\textwidth]{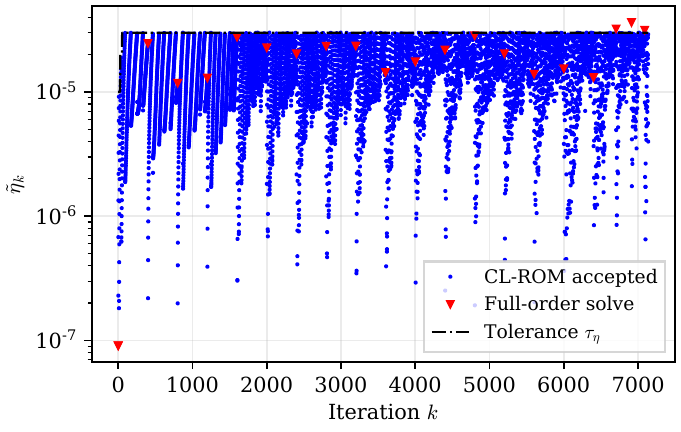}
\caption{Left: objective histories of the full-order baseline and CL-ROM at $h=2\times10^{-3}$. Right: calibrated CL-ROM indicator $\tilde\eta_k=\max(\tilde\eta_s,\tilde\eta_a)$; blue dots mark accepted reduced updates, red triangles full-order updates, and the dash-dotted line the tolerance $\tau_\eta$. The two warm-up full-order updates, at $k=0$ and $1$, have no reduced indicator and are drawn as one triangle at the bottom left.}
\label{fig:hist}
\end{figure}

\begin{table}[htbp]
\centering
\caption{Shape identification: full-order and CL-ROM runs with the settings of Section~\ref{sec:protocol}. The stopping length of the preconditioned full-order run sets the CL-ROM update budget. ``FOM/ROM solves'' counts full-order updates plus verifications / accepted reduced updates; ``corr.'' counts accepted corrections, and ``fact.'' counts Cholesky factorizations. Objectives are averaged over the last fifty updates; timing and counting conventions are given in Section~\ref{sec:protocol}.}
\label{tab:summary}
\small\setlength{\tabcolsep}{4pt}
\begin{tabular}{lrrrrrrrr}
\toprule
method & updates & FOM/ROM solves & corr. & fact. & mean $J$ & time (s) & speedup & $|\Delta J|/J$ \\
\midrule
full order & 7139 & 7139 / -- & -- & 88 & $2.256\times10^{-8}$ & 946.5 & -- & -- \\
CL-ROM & 7139 & $5{+}16$ / 7134 & 821 & 20 & $2.256\times10^{-8}$ & 245.6 & $3.85$ & $2.103\times10^{-6}$ \\
\bottomrule
\end{tabular}
\end{table}

Figure~\ref{fig:phase} shows recovery of the principal target features
from the thin initial heptagon; the narrow trailing tip remains rounded at
finite interface width. The thresholded shape error is $1.568\times10^{-3}$
in area for both full order and CL-ROM. The nearly overlapping objective
histories in Figure~\ref{fig:hist} show that reduced solves preserve the recovery process. In the right panel, the calibrated indicators of the accepted candidates fill the band below the tolerance, which sits at its upper bound for most of the run, and the few full-order solves are scattered along the whole trajectory. Table~\ref{tab:summary} attributes the saving to many
accepted reduced updates and few full-order factorizations.

The full-order run meets the objective test at 7139 updates; the separate
CL-ROM audit run, stopped independently, also stops at 7139 updates. Thus its
stopping length happens to match the budget used for the reported comparison.
Small interfacial oscillations remain, so this criterion describes objective
behavior rather than strict design stationarity. The final CL-ROM phase field differs from the full-order one by
less than $9\times10^{-5}$ in maximum norm, and their thresholded domains
are identical on this grid.


\FloatBarrier
\subsection{Drag minimization in a channel}
\label{sec:flowexp}

The flow problem tests the same decision loop with a nonlinear saddle-point operator. In this setting, an accepted pair may still require an expensive correction, so acceptance alone does not predict runtime. We apply Algorithm~\ref{alg:trrom} to the isolated-body configuration of Test~2 of Gao and Ma \cite{gao2008shape} at $Re = 40$: a solid body of prescribed volume sits in a channel, and its shape is chosen to minimize the drag \eqref{eq:dragflow} exerted on it by the flow, as in the companion Brinkman study \cite{li2026residual}. In \cite{gao2008shape} the objective is instead the total dissipated energy, so the optimal shapes are not expected to coincide. We use the common setting of Section~\ref{sec:protocol}, with one inexact Newton correction after a failed screen, solved by at most eight GMRES iterations preconditioned by the retained LU factors, and a flow-specific POD cutoff of $10^{-12}$. The larger factor-refresh threshold allows the Newton and adjoint solves, which can use about 20 GMRES iterations in one update, to reuse the saved factors.

The channel has height one with parabolic inflow $u = 0.05 - 0.2y^2$, a do-nothing outflow, no-slip walls and $\nu = 0.02/Re$; the diffuse solid volume is constrained to $0.1$ and the initial design is a disk. We pose the problem on the enlarged domain $(-1,3)\times(-\tfrac12,\tfrac12)$, in which the objective is nearly translation invariant, and confine the solid to the interior by a fluid band of width $0.1$ in which the design is fixed. The initial disk uses the signed-distance profile, with its radius adjusted to give the prescribed diffuse solid volume. The sensitivity is not clipped. The prescribed volume is maintained by the exact multiplier \eqref{eq:lambda}, with projection if needed. Stopping follows \eqref{eq:stoptest}. We use $h = 0.025$ and $\varepsilon=0.1$ ($n = 19{,}400$), the penalty constant $C_p=64$, and the pressure regularization $\delta_p=10^{-8}/h$. The sparse LU factorization accepts the small diagonal pivots created by the pressure regularization; rejecting them, as the default pivoting does, causes excessive fill. The resulting factors are slightly inexact, and iterative refinement removes this error. This keeps the full-order baseline efficient.

\begin{figure}[htbp]
\centering
\includegraphics[width=0.47\textwidth]{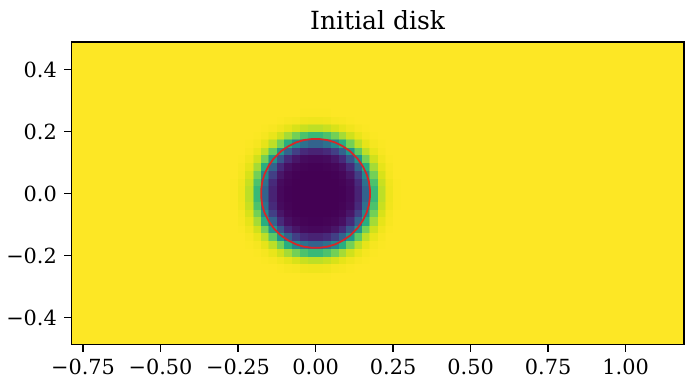}
\hfill
\includegraphics[width=0.47\textwidth]{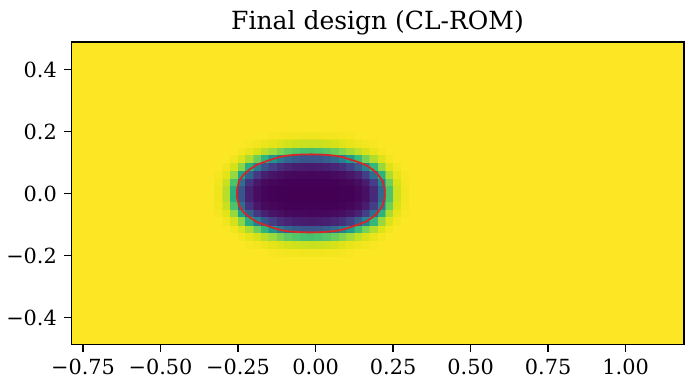}\\[1mm]
\includegraphics[width=0.47\textwidth]{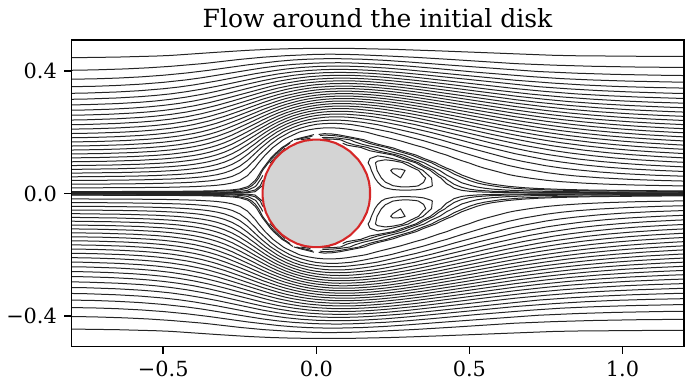}
\hfill
\includegraphics[width=0.47\textwidth]{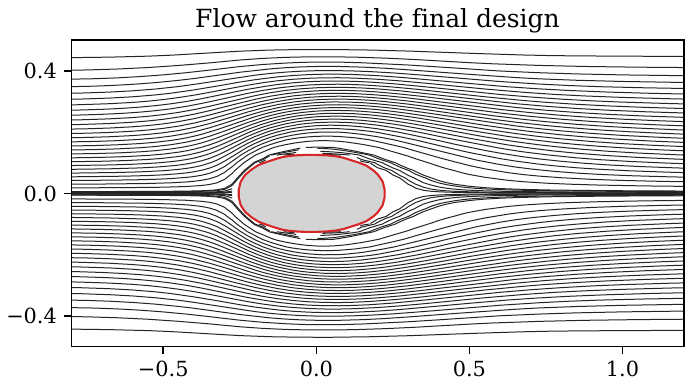}
\caption{Drag minimization at $Re = 40$, $h = 0.025$, $\varepsilon = 4h$, cropped to the body and its near wake, with the $\phi = \tfrac12$ contour in red. Left: the initial diffuse disk and its flow. Right: the CL-ROM design after 5439 updates and its flow. Streamlines are reconstructed from the weighted flux $\phi\bm{u}$, whose divergence is zero up to the pressure regularization. Additional contour levels near the obstacle show how the flow passes around its boundary.}
\label{fig:ddmbody}
\end{figure}

Figure~\ref{fig:ddmbody} shows the initial disk stretching into a symmetric lens as drag decreases. The full-order run satisfies the objective test after 5439 updates and preserves volume to round-off. Its final drag is about $31\%$ below that of the initial disk. Table~\ref{tab:ddmrom} gives the resulting objective and computational work. The fixed-$\varepsilon$ refinement in Appendix~\ref{app:ddm_refinement} checks spatial convergence of the implemented flow model at this interface width.

\begin{table}[htbp]
\centering
\caption{Drag minimization, in the format of Table~\ref{tab:summary}. The full-order run meets its objective test after 5439 updates, and the CL-ROM comparison is run for this same update budget. ``Fact.'' counts LU factorizations of Newton Jacobians.}
\label{tab:ddmrom}
\small\setlength{\tabcolsep}{4pt}
\begin{tabular}{lrrrrrrrr}
\toprule
method & updates & FOM/ROM solves & corr. & fact. & mean $J$ & time (s) & speedup & $|\Delta J|/J$ \\
\midrule
full order & 5439 & 5439 / -- & -- & 12 & $9.39483\times10^{-4}$ & 565.3 & -- & -- \\
CL-ROM & 5439 & $2{+}13$ / 5437 & 54 & 12 & $9.39484\times10^{-4}$ & 73.9 & $7.65$ & $1.31\times10^{-6}$ \\
\bottomrule
\end{tabular}
\end{table}

\begin{figure}[htbp]
\centering
\includegraphics[width=0.49\textwidth]{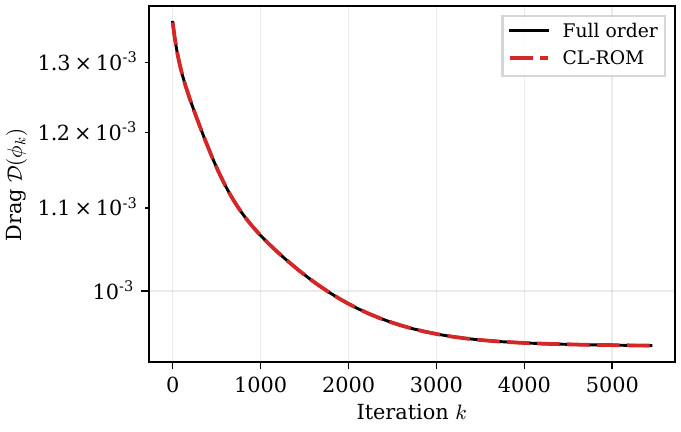}\hfill
\includegraphics[width=0.49\textwidth]{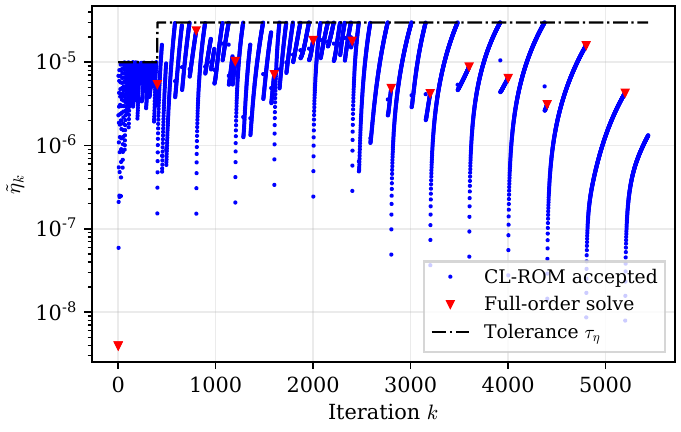}
\caption{The body at $Re = 40$, $h = 0.025$, $\Delta t = 0.01$, in the format of Figure~\ref{fig:hist}. Left: drag histories of the preconditioned full-order baseline and CL-ROM. Right: calibrated indicator $\tilde\eta_k = \max(\tilde\eta_s,\tilde\eta_a)$ of CL-ROM; blue dots mark accepted reduced iterations, red triangles the iterations solved at full order, and the dash-dotted line the tolerance $\tau_\eta$.}
\label{fig:ddmhist}
\end{figure}

Figure~\ref{fig:ddmhist} shows mostly accepted reduced iterations, with occasional full-order refreshes as the body stretches into its final shape. The right panel has the sawtooth pattern of the screen: after each refresh the indicator drops by about two orders of magnitude and then grows as the design moves away from the refreshed basis, and the refreshes become sparser as the design settles.
Table~\ref{tab:ddmrom} shows that CL-ROM tracks the full-order drag
and replaces almost all full-order state--adjoint solves. Only 54 accepted
pairs need a correction on this grid, so most updates use an uncorrected
reduced solution. The resulting speedup is $7.65$.
The full-order and budget-matched CL-ROM final phase fields differ by less
than $6.1\times10^{-6}$ in maximum norm and have the same thresholded
domain. The separate CL-ROM trajectory used for the accepted-pair audit in
Table~\ref{tab:errstudy} applies its own stopping test and stops after 5428
updates, compared with the 5439-update budget of the reported comparison.

\subsection{Compliance minimization of a cantilever}
\label{sec:cantexp}

The cantilever tests whether the same CL-ROM decisions transfer to a traction-free elastic boundary and a self-adjoint objective. We use the two-dimensional cantilever of Takezawa, Nishiwaki and Kitamura \cite[Section 5.1]{takezawa2010shape}. The free structural boundary is traction free through the $\phi$-weighted elasticity model \eqref{eq:elast_ddm}, without a boundary penalty. The rectangular domain $D=(0,2)\times(0,1)$ is clamped on the left and carries a unit downward point load at the center of the right edge. The objective is compliance plus $f_V=80$ times material volume. We use a $200\times100$ mesh ($h=0.01$), plane stress with Young's modulus $E=1$ and Poisson ratio $\nu_{\rm P}=0.3$, so that $\mu=E/(2(1+\nu_{\rm P}))$ and $\lambda=E\nu_{\rm P}/(1-\nu_{\rm P}^2)$. We set $k_{\min}=10^{-5}$ and use SIMP exponent $p=3$. Other CL-ROM parameters follow Section~\ref{sec:protocol}; rejected reduced solves use the same retained-factor CG fallback as the full-order baseline.

\begin{figure}[htbp]
\centering
\includegraphics[width=0.47\textwidth]{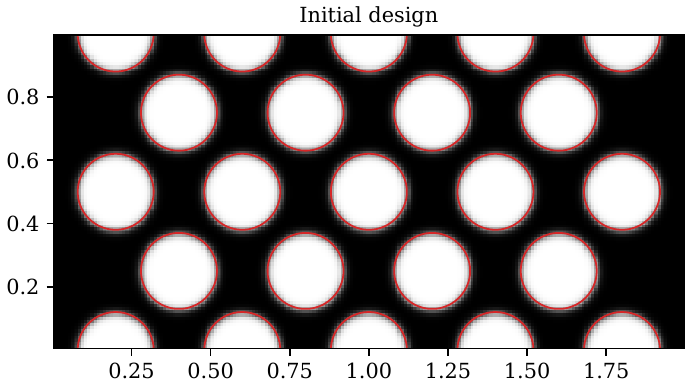}
\hfill
\includegraphics[width=0.47\textwidth]{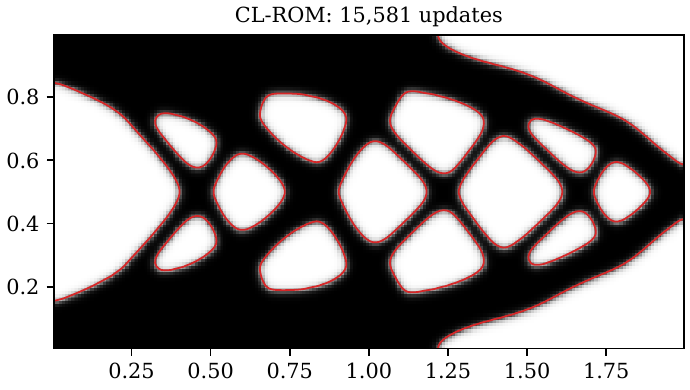}
\caption{The cantilever of \cite[Section 5.1]{takezawa2010shape} in the diffuse-domain setting, $200\times 100$ elements, with the $\phi = \tfrac12$ contour in red: the perforated-plate initial design (left) and the CL-ROM design after the full-order objective test is met at 15581 updates (right). The full-order and CL-ROM thresholded domains are identical on this grid, with a maximum phase-field difference of $2.7\times10^{-3}$.}
\label{fig:cant}
\end{figure}

We initialize the design with the perforated plate of \cite[Fig.~6(b)]{takezawa2010shape}: staggered rows of circular holes leave $59\%$ material (Figure~\ref{fig:cant}, left). This choice matters because the gated design flow does not nucleate holes from a uniform solid: $G'(\phi)=\phi^2(1-\phi)^2$ vanishes at $\phi\equiv1$. The objective test was met at 15581 updates, and CL-ROM is compared at this same total budget.

\begin{figure}[htbp]
\centering
\includegraphics[width=0.49\textwidth]{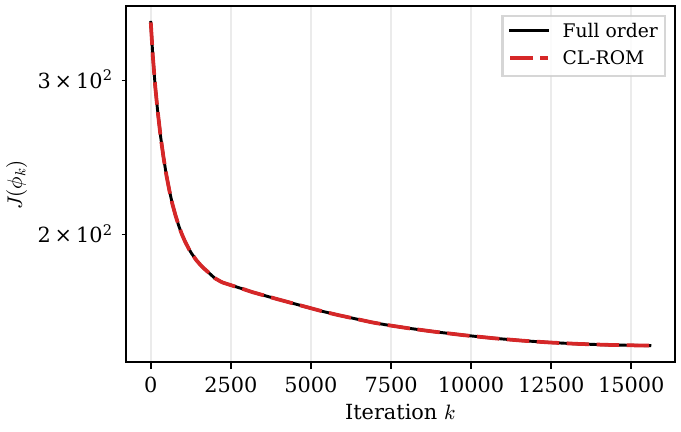}\hfill
\includegraphics[width=0.49\textwidth]{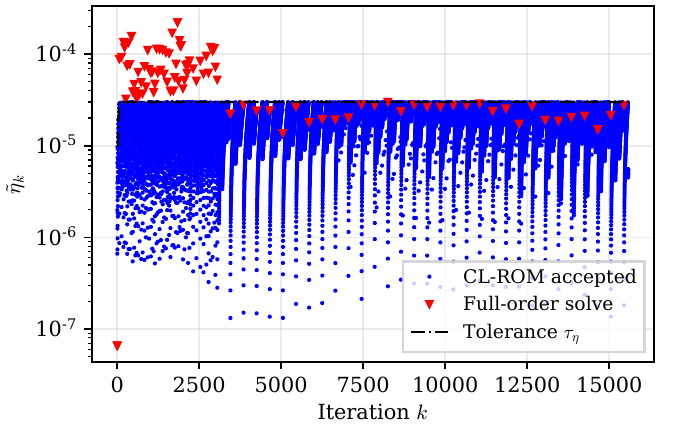}
\caption{The cantilever at matched trajectory length, in the format of Figure~\ref{fig:hist}. Left: objective histories of the preconditioned full-order baseline and CL-ROM. Right: calibrated indicator $\tilde\eta_k$ of CL-ROM; blue dots mark accepted reduced iterations, red triangles the iterations solved at full order, and the dash-dotted line the tolerance $\tau_\eta$.}
\label{fig:canthist}
\end{figure}

\begin{table}[htbp]
\centering
\caption{Compliance minimization of the cantilever, in the format of Table~\ref{tab:summary}. The full-order run meets the objective test at 15581 updates; CL-ROM uses this same budget. The full-order timing includes one restart factorization. ``fact.'' counts stiffness-matrix Cholesky factorizations.}
\label{tab:cantconv}
\small\setlength{\tabcolsep}{4pt}
\begin{tabular}{lrrrrrrrr}
\toprule
method & updates & FOM/ROM solves & corr. & fact. & mean $J$ & time (s) & speedup & $|\Delta J|/J$ \\
\midrule
full order & 15581 & 15581 / -- & -- & 252 & $149.50225$ & 1576.6 & -- & -- \\
CL-ROM & 15581 & $63{+}31$ / 15518 & 258 & 93 & $149.50151$ & 236.6 & $6.66$ & $4.969\times10^{-6}$ \\
\bottomrule
\end{tabular}
\end{table}

Starting from the perforated plate, the design evolves into a truss-like
structure (Figure~\ref{fig:cant}, right): two chords run from the clamped
edge to the loaded point and are joined by diagonal members, and the circular
holes become polygonal voids that shrink toward the tip. The objective histories of the two solvers remain close (Figure~\ref{fig:canthist}). Its right panel shows two phases: during the first few thousand updates, while the holes reshape, candidates often exceed the tolerance and are replaced by full-order solves; afterwards nearly every candidate is accepted with its indicator just below $\tau_\eta$. CL-ROM is faster than full order at the matched
budget (Table~\ref{tab:cantconv}). The final objective evaluated at full order is about $57\%$ below the initial value, and the two thresholded
terminal designs agree. Thus the compact reduced solver preserves the
objective and thresholded design at the full-order stopping horizon.

\FloatBarrier
\subsection{Accuracy of accepted reduced solutions}
\label{sec:audit}

Figures~\ref{fig:hist}, \ref{fig:ddmhist}, and \ref{fig:canthist} show only
that the full-order and reduced trajectories agree. Here we instead measure
the error of the accepted reduced solutions directly: separate
runs re-solve the state and adjoint at full order at every fiftieth design
whose reduced pair passes the residual screen. Samples are taken before any
objective-triggered correction or verification, so they include pairs that
were later replaced. These full-order re-solves do not alter the trajectory and are
excluded from the reported speedups.

\begin{table}[H]
\centering\small
\caption{Same-design accuracy of residual-screened candidates: every fiftieth
of the 7134 and 5426 accepted candidates in the independently stopped
shape-identification and drag audit runs, and of the 15518 accepted
candidates in the matched 15581-update compliance run. State and objective errors are
maximum sampled relative errors; drag-minimization state errors use velocity
components. The last column is the fraction of samples in which the
calibrated state indicator $\tilde\eta_s$ is smaller than the relative error
of the full state vector. For shape identification and drag minimization, these samples come from
separate CL-ROM trajectories stopped by their own tests. Compliance samples
come from the matched 15581-update trajectory.}
\label{tab:errstudy}
\begin{tabular}{lrrrr}
\toprule
problem & samples & state error & objective error & $\tilde\eta_s<$ state error\\
\midrule
shape identification & 142 & $4.4\times10^{-3}$ & $4.4\times10^{-5}$ & $73\%$\\
drag minimization & 108 & $2.5\times10^{-5}$ & $5.0\times10^{-5}$ & $30\%$\\
compliance minimization & 310 & $9.1\times10^{-4}$ & $1.0\times10^{-6}$ & $70\%$\\
\bottomrule
\end{tabular}
\end{table}

Table~\ref{tab:errstudy} confirms that the residual screen is an error
indicator, not an error bound: the calibrated indicator frequently
underestimates the state error, and the sampled state errors can exceed the
screening threshold $\tau_\eta$ by orders of magnitude. The objective error
nevertheless stays below $\tau_J$ in every sample. These results were
obtained with the objective check in Algorithm~\ref{alg:trrom} enabled.
During the optimization, the check triggers corrections and full-order
verifications when the objective estimate \eqref{eq:dwr} is large and
adjusts $\tau_\eta$ from the observed objective error.

\FloatBarrier
\section{Discussion and Conclusions}
\label{sec:discussion}

This work extends CL-ROM from Brinkman-penalized flow \cite{li2026residual} to diffuse-domain shape optimization. The main contribution is to accommodate the different operators, diffuse boundary terms and sensitivities of the chosen DDM formulations within a common reduction and verification procedure. The fixed background grid allows snapshots from changing geometries to be combined directly. Shape identification, nonlinear drag minimization and self-adjoint compliance minimization demonstrate this extension for problems with different boundary conditions and state--adjoint structures.

At matched numbers of design updates, CL-ROM closely follows the full-order objective histories and final designs, with speedups of $3.85$--$7.65$ over a full-order baseline that reuses factorizations as preconditioners. These gains come from replacing most full-order solves by reduced solves, supplemented by corrections and occasional full-order verification. For the symmetric problems, Proposition~\ref{prop:ddm_local} bounds the reduced-state error in terms of changes in the operator and load, snapshot error, and POD truncation. It clarifies the conditions under which recent, corrected snapshots provide an effective local space; the nonlinear flow case lies outside its scope. The accuracy checks also identify the limits of the approach. Calibrated residual indicators can underestimate the state error, so residual acceptance alone does not certify solution accuracy. With objective checks enabled, all sampled objective errors remain below the prescribed tolerance. These observations support the combined screening, correction and verification strategy, without providing an error guarantee at every update.

Further speedup is limited by work that still depends on the full grid: reduced assembly, residual evaluation, sensitivity evaluation and the design update. Hyperreduction of selected operator and load terms could lower assembly and residual-evaluation costs. Its additional approximation error would need to be assessed in the screening procedure and monitored through full-residual checks. Complementary improvements to sensitivity evaluation and the design update could further reduce the total optimization cost.

\FloatBarrier
\section*{Acknowledgments}
Z.W.\ was partially supported by the U.S. National Science Foundation under grant numbers DMS-2038080, DMS-2245097, and DMS-2513112. L.J.\ was partially supported by the U.S. National Science Foundation under grant number DMS-2409634.

\FloatBarrier

\appendix
\makeatletter
\renewcommand{\@seccntformat}[1]{%
  \ifstrequal{#1}{section}{Appendix~\thesection.\quad}{\csname the#1\endcsname\quad}}
\makeatother
\FloatBarrier
\section{Derivation for the Shape-Identification Case}
\label{app:derivation}

The state, adjoint and sensitivity follow from one Lagrangian construction.
Let $V=H^1_0(D)$ and $e=\phi u-\bar u$.

\paragraph{State equation and Lagrangian.}
For a test function $\zeta\in V$, define the weak state residual
\begin{equation*}
R(u,\phi;\zeta)=\int_D
k(\phi)\nabla u\cdot\nabla\zeta
+\lambda_\varepsilon\delta_\sigma(\phi)u\zeta
-\phi f\zeta\,\mathrm{d}\bm{x}.
\end{equation*}
The Lagrangian is
\[
L(u,w;\phi)=\int_D e^2\,\mathrm{d}\bm{x}-R(u,\phi;w),
\]
where the adjoint $w\in V$ is the Lagrange multiplier. Variation with
respect to $w$ in a direction $\zeta$ gives
$\partial_wL[\zeta]=-R(u,\phi;\zeta)$; stationarity therefore returns the
weak form of the state equation.

\paragraph{Adjoint equation.}
Variation with respect to $u$ in a direction $v\in V$ gives
\[
\partial_uL[v]=\int_D
2\phi e\,v-k(\phi)\nabla v\cdot\nabla w
-\lambda_\varepsilon\delta_\sigma(\phi)vw\,\mathrm{d}\bm{x}.
\]
Setting this variation to zero for every $v$ yields
\eqref{eq:adjoint_cont}; the same symmetric operator appears in the state and
adjoint equations, with a different right-hand side.

\paragraph{Design sensitivity.}
Once the state and adjoint equations hold, their variations cancel from the
total derivative. For a design direction $\vartheta$, the remaining explicit
variation uses
\[
\delta_\sigma'(\phi)[\vartheta]
=\frac{\nabla\phi\cdot\nabla\vartheta}
{\sqrt{|\nabla\phi|^2+\sigma^2}}
\]
and is
\[
\partial_\phi L[\vartheta]=\int_D
\big[2ue+fw-(1-k_{\min})\nabla u\cdot\nabla w\big]\vartheta
-\lambda_\varepsilon uw\,
\frac{\nabla\phi\cdot\nabla\vartheta}
{\sqrt{|\nabla\phi|^2+\sigma^2}}\,\mathrm{d}\bm{x}.
\]
This is the sensitivity \eqref{eq:sens_cont}.

\FloatBarrier
\section{Derivation for the Drag-Minimization Case}
\label{app:flow}

The derivation uses the boundary-flux objective \eqref{eq:dragflow} and
the weighted state equations \eqref{eq:nse_ddm}. We assume sufficient
regularity for the following integrations by parts. The phase field is
fixed at one in a neighborhood of $\partial D$, and admissible design
variations vanish there.

\paragraph{State equation and Lagrangian.}
Write the momentum and continuity residuals as
\[
\bm R_m=\phi(-\nu\Delta\bm u+\nabla p)
 +(\bm u\cdot\nabla)(\phi\bm u)+a\bm u,
\qquad R_c=\nabla\cdot(\phi\bm u),
\qquad a=\alpha_\varepsilon(1-\phi)^2.
\]
On fields satisfying the prescribed state boundary conditions, introduce
\[
L(\bm u,p,\bm w,\pi;\phi)
=J(\bm u,p)-\int_D(\bm w\cdot\bm R_m+\pi R_c)\,\mathrm d\bm x.
\]
Variations in the multipliers give the state equations. The drag contains
boundary pressure and normal derivatives, so deriving its adjoint requires
retaining the boundary terms in the state variation.

\paragraph{Adjoint equation.}
A state variation $(\bm v,q)$ satisfies $\bm v=\bm0$ on
$\Gamma_i\cup\Gamma_w$ and $\nu\partial_{\bm n}\bm v-q\bm n=\bm0$
on $\Gamma_o$. Integrating the variation of the residual term by parts gives
\[
\begin{aligned}
\delta\!\int_D(\bm w\cdot\bm R_m+\pi R_c)\,\mathrm d\bm x
={}&\int_D\bm v\cdot\Big[-\nu\Delta(\phi\bm w)
 +[\nabla(\phi\bm u)]^{\!\top}\bm w
 -\phi\big((\bm u\cdot\nabla)\bm w+(\nabla\cdot\bm u)\bm w\big)
 +a\bm w-\phi\nabla\pi\Big]\,\mathrm d\bm x\\
&-\int_Dq\,\nabla\cdot(\phi\bm w)\,\mathrm d\bm x\\
&+\int_{\partial D}\Big[
 \big(\nu\partial_{\bm n}(\phi\bm w)+\phi(\bm u\cdot\bm n)\bm w
      +\phi\pi\bm n\big)\cdot\bm v
 -\nu\phi\bm w\cdot\partial_{\bm n}\bm v
 +\phi(\bm w\cdot\bm n)q\Big]\,\mathrm ds.
\end{aligned}
\]
The two convection contributions come from varying the advecting and
advected velocities. The objective variation, after using the state
boundary conditions, is
\[
\delta J
=\int_{\Gamma_i}\big(q+\nu\partial_{\bm n}v_x\big)\,\mathrm ds
 +\nu\int_{\Gamma_w}\partial_{\bm n}v_x\,\mathrm ds
 -2\int_{\Gamma_o}u_xv_x\,\mathrm ds.
\]
At the inlet, $\partial_x=-\partial_{\bm n}$ and $\bm v=\bm0$.
At the outlet, $q-\nu\partial_xv_x=0$, so the pressure and normal
viscous variations cancel. Setting $\bm w=-\bm e_x$ on
$\Gamma_i\cup\Gamma_w$ matches the inlet pressure and the Dirichlet-boundary
normal-derivative terms. At $\Gamma_o$, the last two residual boundary
terms cancel by the linearized state traction condition. Equating the
remaining coefficient of $\bm v$ to $-2u_x\bm e_x$, and setting the
interior coefficients of $\bm v$ and $q$ to zero, yields
\eqref{eq:adjoint_flow}. In particular, imposing homogeneous adjoint
velocity on the inlet and walls would omit part of the drag variation.

\paragraph{Design sensitivity.}
At a stationary state--adjoint pair, differentiating $L$ in a design
direction $\vartheta$ leaves only the explicit coefficient dependence:
\[
\begin{aligned}
\partial_\phi L[\vartheta]=-\int_D\Big[&
 \nu\nabla\bm u:\nabla(\vartheta\bm w)
 +\big((\bm u\cdot\nabla)(\vartheta\bm u)\big)\cdot\bm w
 -2\alpha_\varepsilon(1-\phi)\vartheta\,\bm u\cdot\bm w\\
 &-p\nabla\cdot(\vartheta\bm w)
 +\pi\nabla\cdot(\vartheta\bm u)\Big]\,\mathrm d\bm x.
\end{aligned}
\]
Here integration by parts in the viscous and pressure terms produces no
boundary contribution because $\vartheta$ vanishes near $\partial D$.
This gives \eqref{eq:sens_flow}. Adding the fluid-volume constraint
adds a constant volume multiplier to the sensitivity; its gated discrete
counterpart is determined by \eqref{eq:lambda}.

\FloatBarrier
\section{Derivation for the Compliance-Minimization Case}
\label{app:elast}

The same Lagrangian construction gives the elastic adjoint and sensitivity.
The elasticity operator is symmetric, so compliance is self-adjoint.

\paragraph{State equation and Lagrangian.}
Let $V = \{\bm{v} \in [H^1(D)]^d : \bm{v} = \bm{0} \text{ on } \Gamma_D\}$.
Multiplying \eqref{eq:elast_ddm} by $\bm{v} \in V$ and integrating the principal
part by parts gives
\[
\int_D \bm{v}\cdot\Big(-\nabla\cdot\big(k(\phi)\,A\,\bm{S}(\bm{u})\big)\Big)\,\mathrm{d}\bm{x}
= \int_D k(\phi)\,A\,\bm{S}(\bm{u}) : \nabla\bm{v}\,\mathrm{d}\bm{x}
 - \int_{\partial D} k(\phi)\,\big(A\,\bm{S}(\bm{u})\big)\bm{n}\cdot\bm{v}\,\mathrm{d}s .
\]
Because $A\bm{S}(\bm{u})$ is symmetric, $A\bm{S}(\bm{u}):\nabla\bm{v} =
A\bm{S}(\bm{u}):\bm{S}(\bm{v})$, which is what makes the resulting bilinear form
symmetric. The clamp is built into $V$, and the fixed traction supplies the
boundary integral on $\Gamma_N$. In the sharp-interface limit with vanishing
void stiffness, the design boundary carries no stress. These conditions give
\eqref{eq:elast_weak_design}. Define
\begin{equation*}
a_\phi(\bm{u},\bm{v}) = \int_D k(\phi)\,A\,\bm{S}(\bm{u}) : \bm{S}(\bm{v})\,\mathrm{d}\bm{x},
\qquad
\ell(\bm{v}) = \int_D \phi\,\bm{f}\cdot\bm{v}\,\mathrm{d}\bm{x} + \int_{\Gamma_N}\bm{t}\cdot\bm{v}\,\mathrm{d}s,
\end{equation*}
so that the state problem reads $a_\phi(\bm{u},\bm{v}) = \ell(\bm{v})$ for all
$\bm{v} \in V$, with $a_\phi$ symmetric and, under the floor $k_{\min}$ of
Section~\ref{sec:elast}, coercive on $V$ when the clamp has positive measure.
For an objective $J(\bm{u},\phi)$, introduce the multiplier $\bm{w}\in V$
and the Lagrangian
\[
L(\bm{u},\bm{w};\phi)
=J(\bm{u},\phi)-\big(a_\phi(\bm{u},\bm{w})-\ell(\bm{w})\big).
\]
Variation with respect to $\bm{w}$ returns the weak state equation.

\paragraph{Adjoint equation.}
For compliance,
$J_{\rm c}(\phi) = \ell(\bm{u})$ by \eqref{eq:compliance}. The state variation in
a direction $\bm{v} \in V$ is
\[
\partial_{\bm{u}} L[\bm{v}] = \ell(\bm{v}) - a_\phi(\bm{v},\bm{w}),
\]
so stationarity requires $a_\phi(\bm{w},\bm{v}) = \ell(\bm{v})$ for all $\bm{v}$,
using the symmetry of $a_\phi$. This is the state problem itself: $\bm{w} =
\bm{u}$, and no second solve is required. Thus the reduced model in Section~\ref{sec:trrom} requires only one basis
for compliance, whereas the other two problems use separate state and adjoint bases.

\paragraph{Design sensitivity.}
At a stationary state--adjoint pair, differentiating $L$ with respect to the
design in a direction $\vartheta$ leaves only the explicit $\phi$-dependence. Using the
interpolated stiffness $k(\phi)A$ from \eqref{eq:elast_weak_design},
\[
\partial_\phi L[\vartheta]
= \underbrace{2\int_D \vartheta\,\bm{f}\cdot\bm{u}\,\mathrm{d}\bm{x}}_{\text{from the objective and state load}}
 - \int_D \vartheta\,k'(\phi)\,A\,\bm{S}(\bm{u}) : \bm{S}(\bm{u})\,\mathrm{d}\bm{x} ,
\]
and if the body force vanishes or acts only where the design is fixed, the
first term vanishes and \eqref{eq:sens_compliance} follows. This is the case
in the experiments of Section~\ref{sec:cantexp}, where the load is a
boundary traction. The sign is worth noting: since $k' > 0$ and
$A\bm{S}(\bm{u}):\bm{S}(\bm{u}) \geq 0$ is twice the local strain energy density,
$J_{\rm c}' \leq 0$ everywhere. Adding material can only reduce compliance, so
the volume term in \eqref{eq:compliance} penalizes excess material.

The volume term in \eqref{eq:compliance} adds
$f_V\int_D\vartheta\,\mathrm d\bm x$ to this variation.

\FloatBarrier
\FloatBarrier
\section{Verification of the Three Diffuse Models}
\label{app:ddm_refinement}

We check spatial refinement of each implemented diffuse model on a fixed
two-dimensional geometry. The interface width is held at the value used in
the corresponding optimization: $\varepsilon_0=0.008$, $0.1$, and $0.04$
for shape identification, drag minimization, and compliance minimization, respectively. With
$h_0=0.002$, $0.025$, and $0.01$, we solve at
$h=2h_0,h_0,h_0/2,h_0/4$ and compare each coarser solution with the
$h_0/4$ solution of the \emph{same} diffuse problem. The circular interfaces
use the signed-distance profile \eqref{eq:phase}; $\phi\approx1$ inside
the disk for shape identification and compliance minimization and outside the obstacle for drag minimization.

For the shape-identification model, $\Omega$ is the disk $r<R=0.075$ in
$D=(-0.12,0.12)^2$. We use $f=4$ in $\Omega$ and zero Dirichlet data on
$\Gamma$, for which the sharp solution is $u_\star=(R^2-r^2)_+$.
For compliance minimization, $R=0.25$, $D=(-0.4,0.4)^2$, and the body force in
$\Omega$ is derived from the sharp field
$\bm u_\star=q(r^2)(x,y)$, where $q(s)=(1-s/R^2)_+^3$.
Both this field and its stress vanish on $\Gamma$; the outer box is fixed.
We use the plane-stress tensor, SIMP exponent, and void floor of the
optimization model.
For drag minimization, $\Omega=D\setminus\overline{B_R}$ with
$D=(-0.75,0.75)^2$ and $R=0.25$. We set $p_\star=0$ and
$\bm u_\star=g(r^2)(-y,x)$, where
$g(s)=(s-R^2)(B^2-s)^3$ for $R^2<s<B^2$, zero otherwise,
and $B=0.6$. This divergence-free field vanishes on the body and the
outer box and is tangent to the circular phase contours. The forcing
$\bm f=-\nu\Delta\bm u_\star+
(\bm u_\star\cdot\nabla)\bm u_\star$ is prescribed in the fluid, with
$\phi\bm f$ added to the diffuse momentum equation
\eqref{eq:nse_ddm}. The nonlinear MAC operator matches the flow experiment. To isolate
mesh refinement, the pressure regularization is held at its main-grid
value $\delta_p=4\times10^{-7}$ on all four grids.

Note that, in each case, the sharp solution serves only to define the source term; since
at fixed $\varepsilon_0$ the diffuse solution converges to a diffuse limit
that differs from the sharp solution by $O(\varepsilon_0)$, the errors below
are measured against the $h_0/4$ solution of the same diffuse problem.

Table~\ref{tab:ddm_refinement} reports relative discrete
$L^2(\Omega)$ state differences. These differences
decrease at approximately second order as $h$ is refined at fixed
$\varepsilon_0$. The last measured rates also reflect the finite resolution
of the reference. These tests support the spatial discretizations at the
selected interface widths; they do not measure sharp-interface error or
establish grid convergence of an optimized shape.

\begin{table}[htbp]
\centering\small
\caption{Two-dimensional mesh refinement at each optimization interface
width $\varepsilon_0$. Errors are measured against the solution at $h_0/4$ with the
same $\varepsilon_0$. The last column gives
$\log_2(E_{2h_0}/E_{h_0})$ and $\log_2(E_{h_0}/E_{h_0/2})$.}
\label{tab:ddm_refinement}
\begin{tabular}{lrrrrrc}
\toprule
problem & $\varepsilon_0$ & $h_0$ & $E_{2h_0}$ & $E_{h_0}$ & $E_{h_0/2}$ & rates \\
\midrule
shape identification & $0.008$ & $0.002$ & $2.41\times10^{-2}$ & $6.03\times10^{-3}$ & $1.20\times10^{-3}$ & $2.00/2.33$ \\
drag minimization & $0.1$ & $0.025$ & $4.29\times10^{-2}$ & $8.72\times10^{-3}$ & $2.10\times10^{-3}$ & $2.30/2.05$ \\
compliance minimization & $0.04$ & $0.01$ & $5.98\times10^{-3}$ & $1.65\times10^{-3}$ & $3.27\times10^{-4}$ & $1.85/2.34$ \\
\bottomrule
\end{tabular}
\end{table}

\FloatBarrier
\section{Sensitivity to the Design Time Step}
\label{sec:tau}

Table~\ref{tab:tau_common} compares $\Delta t=0.01$ and $\Delta t=0.05$,
keeping the initial designs, calibrated coefficients,
grids, stopping rule, and state solvers fixed. CL-ROM accelerates all three
problems at both step sizes. The larger step substantially reduces the
optimization time for the flow and the cantilever while giving comparable objective
levels. For shape identification, however, it worsens the recovered shape
and gives little reduction in CL-ROM time. These results support the
conservative common choice $\Delta t=0.01$ in the main comparisons and show
that the acceleration persists when the full-order design iteration uses
a larger step.

\begin{table}[htbp]
\centering
\caption{Step-size sensitivity; all other settings as in Section~\ref{sec:protocol}. The drag-minimization rows use $h=0.025$. Mean $J$ is the full-order mean over the last fifty updates; $|\Delta J|/J$ is the relative difference between the CL-ROM and full-order means, and ``band'' is the final CL-ROM 10--90\% transition width. The $\Delta t=0.01$ compliance-minimization FOM time includes one restart factorization.}
\label{tab:tau_common}
\small\setlength{\tabcolsep}{4pt}
\begin{tabular}{lcrrrrrrr}
\toprule
problem & $\Delta t$ & updates & FOM (s) & CL-ROM (s) & speedup & mean $J$ & $|\Delta J|/J$ & band \\
\midrule
shape identification & $0.01$ & 7139 & 946.5 & 245.6 & $3.85$ & $2.2561\times10^{-8}$ & $2.1\times10^{-6}$ & $2.67h$ \\
 & $0.05$ & 5324 & 793.4 & 241.6 & $3.28$ & $4.2726\times10^{-8}$ & $9.9\times10^{-6}$ & $2.69h$ \\
\midrule
drag minimization & $0.01$ & 5439 & 565.3 & 73.9 & $7.65$ & $9.39483\times10^{-4}$ & $1.3\times10^{-6}$ & $2.62h$ \\
in a channel & $0.05$ & 2029 & 194.9 & 30.1 & $6.47$ & $9.35542\times10^{-4}$ & $1.5\times10^{-5}$ & $2.66h$ \\
\midrule
compliance minimization & $0.01$ & 15581 & 1576.6 & 236.6 & $6.66$ & $149.5023$ & $5.0\times10^{-6}$ & $2.75h$ \\
of a cantilever & $0.05$ & 3848 & 477.1 & 94.4 & $5.05$ & $149.2166$ & $2.1\times10^{-6}$ & $2.83h$ \\
\bottomrule
\end{tabular}
\end{table}

\bibliographystyle{plainnat}
\bibliography{ref_ddm}

\end{document}